\documentclass{amsart}
\usepackage{amssymb}
\usepackage{xspace}
\usepackage{amsmath}
\usepackage{relsize}
\usepackage{nicefrac} 
\usepackage{marginnote}
\usepackage{enumitem}
\usepackage[pdftex]{graphicx}
\usepackage{epstopdf}
\usepackage{fancyvrb}
\usepackage{subfig}
\usepackage{float}
\usepackage{multirow}
\usepackage{booktabs}
\usepackage[square]{natbib}
\usepackage{hyperref}
\usepackage{tikz}

\tikzset{
    declare function={ceiling(\x)=round(\x+0.5);}
}

\AppendGraphicsExtensions{.tif}

\newcommand{\e}{\varepsilon}
\newcommand{\R}{\ensuremath \mathbb R}

\newcommand{\Z}{\ensuremath{\mathbb Z}}
\renewcommand{\P}{\ensuremath P}

\newcommand{\Ncal}{\ensuremath \mathcal N}

\newcommand{\oas}{\ensuremath{o_{a.s.}}}
\newcommand{\OP}{\ensuremath{O_{P}}}
\newcommand{\oP}{\ensuremath{o_{P}}}

\newcommand{\indicator}[1]{\ensuremath \mathbf{1}_{\{#1\}}}

\newcommand{\X}{\ensuremath{\mathbb X}}

\newcommand{\Xone}{\ensuremath{ X^{(1)}}}
\newcommand{\Xtwo}{\ensuremath{ X^{(2)}}}
\newcommand{\Xz}{\ensuremath{ X^{(z)}}}

\newcommand{\V}{\ensuremath{\mathbb V}}
\newcommand{\Y}{\ensuremath{\mathbb Y}}
\newcommand{\Yone}{\ensuremath{ \Y^{(1)}}}
\newcommand{\Ytwo}{\ensuremath{ \Y^{(2)}}}
\newcommand{\Yz}{\ensuremath{ \Y^{(z)}}}

\newcommand{\cp}{m}

\newcommand{\cpesti}{\hat m}

\newcommand{\D}{\mathcal{D}}

\newcommand{\ParaSpace}{\ensuremath{\Theta}}
\newcommand{\Para}{\ensuremath{\theta}}

\newcommand{\estiPara}{\ensuremath{\hat \Para_N}}

\newcommand{\limPara}{{\ensuremath{\tilde \Para}}}

\newcommand{\lossFct}{\mathcal{E}}

\newcommand{\transpose}[1]{#1^{\mathtt{T}}}

\newcommand{\norm}[1]{\left\| #1 \right\|}
\newcommand{\normSquared}[1]{\left\| #1 \right\|^2}

\newcommand{\normEstimatorof}[2]{{\norm{ #1 }_{#2}}}
\newcommand{\normEstimator}[1]{{\normEstimatorof{#1}{A}}}

\newcommand{\normEstimatorofSquared}[2]{{\norm{ #1 }^2_{#2}}}
\newcommand{\normEstimatorSquared}[1]{{\normEstimatorofSquared{#1}{A}}}

\newcommand{\ovl}[1]{\overline{#1}}
\newcommand{\sm}[1]{\ovl {#1}_n}

\newcommand{\HO}{\ensuremath{H_0}\xspace}
\newcommand{\HI}{\ensuremath{H_1}\xspace}

\newcommand{\E}{\ensuremath{\mathbb E}}

\newcommand{\argmax}{\operatorname*{arg\,max}}

\newcommand{\dconv}{ \ensuremath{ \overset{d}{\longrightarrow}}}
\newcommand{\pconv}{\ensuremath{ \overset p {\longrightarrow}}}
\newcommand{\asconv}{\ensuremath{ \overset {a.s} {\longrightarrow}}}

\newcommand{\HajekRenyi}{H{\'{a}}jek-R{\' e}nyi\xspace}

\newcommand{\teststatistic}{{test statistic}\xspace}

\newcommand{\timeseries}{{time-series}\xspace}

\newcommand{\absatz}{\\[5mm]}

\newtheorem{Definiti}{Definition}[section]

\newtheorem{Remar}{Remark}[section]

\newtheorem{Propo}{Proposition}[section]
\newenvironment{Proposition}{\begin{Propo}}{\end{Propo}}

\newtheorem{Theore}{Theorem}[section]
\newenvironment{Theorem}{\begin{Theore}}{\end{Theore}}

\newtheorem{Lemm}{Lemma}[section]
\newenvironment{Lemma}{\begin{Lemm}}{\end{Lemm}}

\newcommand{\setmylabel}[2]{\textbf{#1.#2}}

\newcounter{saveassumptions}
\newenvironment{MyAssumption}
{ \begin{enumerate}[labelindent=\parindent ,leftmargin=*, label=\setmylabel{A}{{\arabic*}}]\setcounter{enumi}{\value{saveassumptions}} }
{ \setcounter{saveassumptions}{\value{enumi}} \end{enumerate} }

\newcounter{savenotations}
\newenvironment{MyNotation}
{ \begin{enumerate}[labelindent=\parindent ,leftmargin=*, label=\setmylabel{N}{{\arabic*}} ]\setcounter{enumi}{\value{savenotations}} }
{ \setcounter{savenotations}{\value{enumi}} \end{enumerate} }

 \newtheorem{Corollary}{Corollary}[section]

\newenvironment{Proposition*}{\textbf{ Proposition }}{}

 \newenvironment{Theorem*}{\textbf{ Theorem }}{}
\newenvironment{Theoremwith*}[1]{\textbf{ Theorem \normalfont (#1)}\begin{itemize}\item[]}{\end{itemize}}

 \newenvironment{Corollary*}{\textbf{ Corollary }}{}

 \newenvironment{Lemma*}{\textbf{ Lemma }}{}

\newenvironment{Proof}[1]{\textbf{ Proof #1:}\\}{\hspace*{\fill}\rule{2mm}{2mm}\\[3mm]}
\newenvironment{Proof*}{\textbf{ Proof:}\\}{\hspace*{\fill}\rule{2mm}{2mm}\\[3mm]}

\newcommand{\ls}{\leq}

\usepackage{amsaddr}
\title[Change points for NL(A)R-processes with Neural Networks]{Estimation of Change Points for Non-linear (auto-)regressive processes using Neural Network Functions}

\author{Claudia Kirch}
\email{claudia.kirch@ovgu.de}
\address{Otto-von-Guericke-University Magdeburg, Germany, Department of Mathematics, Institute of Mathematical Stochastics}

\author{Stefanie Schwaar}
\address{Fraunhofer ITWM, Fraunhofer Platz 1, Kaiserslautern}
\address{\emph{and} University of Applied Sciences (HTW Berlin), Berlin, Germany}
\email{stefanie.schwaar@htw-berlin.de}
\begin{document}
\thispagestyle{empty}

\pagenumbering{arabic}

\date{\today}
\begin{abstract} 
In this paper, we propose a new test for the detection of a change in a non-linear (auto-)regressive time series as well as a  corresponding estimator for the unknown time point of the change. To this end, we consider an at-most-one-change model and approximate the unknown (auto-)regression function by a neuronal network with one hidden layer.
It is shown that the test has asymptotic power one for a wide range of alternatives not restricted to changes in the mean of the time series. Furthermore, we prove that the
corresponding estimator converges to the true change point with the optimal rate $\OP(1/n)$ and derive the asymptotic distribution. Some simulations illustrate the behavior of the estimator with a special focus on the misspecified case, where the true regression function is not given by a neuronal network. Finally, we apply the estimator to some financial data. 
\end{abstract}

\maketitle

\section{Introduction}\label{sec:NLAR-NeuNet_Motivation}
Many practical applications depend on the collection of data over time to be used for modeling and subsequently the estimation of the unknown parameters. The consequent analysis as e.g.\ the pricing of financial products usually depends on the precision and correctness of the models and estimators. On the other hand, it is natural, that changes in the structure occur, whenever data is observed over a longer period of time, rendering the subsequent analysis meaningless. As a result the detection and estimation of structural breaks is of great importance in many areas of science,
including  industrial  quality  management (production  monitoring \cite{PageContinuous54}), finance (stock prices, exchange rates), climate studies (global warming \cite{GallagherEtAl13}),  medicine (on-line-monitoring of intensive-care patients \cite{FriedImhoff04}), neuroscience (\cite{astonfmri,stoehr2020detecting}) or geoscience (annual water volume of rivers \cite{gao2010trend}).

While the literature on change point analysis already dates back to the 1950ies (\cite{PageContinuous54}), there has been an increased interest in the area in the last two decades with the focus shifting to more complicated situations including more complicated underlying data structures, see e.g.\ the recent survey papers \cite{AueHorvath13,aminikhanghahi2017survey,truong20,ChoKirch21}.

Classically, a large part
 of the literature is concerned with linear models, where the first papers used CUSUM statistics based on residuals (\cite{Kulperger1985,Horvath1993}) which are computationally easy and have no size problems but only have power against particular kinds of alternatives, effectively alternatives that go along with a mean change in the observations. 
A second generation of tests are universal in the sense that they do have power against all kinds of changes in the parameters, where \cite{davisautoregressive} used a likelihood-ratio test statistics while \cite{HuskovaEtAl07} used the corresponding (equivalent but easier to compute) score type statistic. 
To move away from the linear structure assumed in these papers
\cite{KirchTadjuidjeTesting} develop a testing procedure in non-linear autoregressive processes based on the sample residuals obtained from approximating the unknown auto-regression function with a neural network. 
Similarly as in the linear case, these tests can only detect certain changes, essentially such changes that go along with a mean change in the observations. Therefore, in this work, 
 we generalize the test statistic from ~\cite{KirchTadjuidjeTesting} to a universal test statistic   that  can detect any changes for which the best-approximating neural network model changes. In particular, if the auto-regression function is indeed given by a neural network, all parameter changes will be detected. 
 Additionally, the theory goes beyond the purely autoregressive setup allowing for a regression vector that can possibly depend on exogenous as well as endogenous regressors.
Most prominently, we derive the consistency of the change point estimator with optimal rate as well as  its asymptotic distribution.

The paper is organized as follows: In the next section, we present the modeling of an AMOC (At Most One Change) process using a neural network, along with the underlying model assumptions. Subsequently, we introduce the change point detection procedure, first discussing the test statistic and then the change point estimator, with a focus on employing a more flexible approach for both. The performance of the test and estimator is evaluated through simulations in a controlled setting and demonstrated on a realistic example. Proofs and additional technical details are provided in the appendix.

\section{Modeling with neural networks}
In this section we first describe the modelling  in detail, before giving some results on the corresponding neural network parameters.

\subsection{Model}\label{sec:NLAR-NeuNet_MODEL}
We are interested in analyzing a non-linear (auto-)regressive model, i.e.\ we allow for autoregressive or exogenous regressive models, or a combination of them. To simplify the notation, we collect both autoregressive $\X_t$ and exogenous regressors $ \V_t$ in one vector $\Y_t$:
\begin{MyNotation}
\item Let $\Y_t=\transpose{(\X_t,\V_t)}\in \R^{p+d}$ be a vector containing in the first $p$ coordinates the autoregression vector ($\X_t=\transpose{(X_{t-1},\dots,X_{t-p})}\in \R^p$) and in the last $d$ coordinates the exogenous regression vector ($\V_t=\transpose{(V_{t-1},\dots,V_{t-d})}\in \R^d$).\label{notation:NLAR-NeuNet_combined_regression_vector}
\end{MyNotation}

Because for the autoregressive part of order $p$, the vector of autoregressors $\X_t$ is only fully available after the first $p$ observations, for simplicity of notation, we renumber the observations to $X_{-p+1},\dots,X_0,X_1,\dots,X_n$, thus having the complete regressors for observations $X_1,\ldots,X_n$, where $n$ is the number of relevant observations.

With this notation, we are now ready to state the (auto-)regressive change point model that will be investigated in this work:
\begin{equation} X_t =\begin{cases}
   \Xone_t=g_1(\Yone_t)+\e_t & t\le \cp\,,\\
    \Xtwo_t=g_{2}(\Ytwo_t)+\e_t & t> \cp\,,\\
  \end{cases}\label{eq:NLAR-NeuNet_started_auto-regressive_model}
\end{equation}
where 
$1< \cp=\cp(n)=\lfloor \lambda n \rfloor\le n$, $0<\lambda\ls 1$ is called the change point and the regression functions $g_1$
and $g_{2}$ are assumed to be unknown.

We use a single-layer neural network as an approximation for the non-linear regression function allowing us to construct change point methods based on least-squares estimation.
The use of a shallow neural network is motivated by two things: First, even one-layer neural networks have a universal  approximation property  (see \cite{Horniketal89}), i.e.\ they can approximate any sufficiently smooth regression function  to any degree of accuracy. The proposed approach essentially transforms the problem into a multivariate mean-change problem (of the transformed time series) with a dimension that is equal to the number of unknown parameters in the model. Thus, secondly, we aim to keep this number relatively small by choosing a single-layer network with a moderate number of hidden neurons to avoid   size problems. Additionally, this may even help with power as long as the change is clearly visible by the best-approximating neural networks, see \cite{kirch2015eeg} for a detailed discussion of this phenomenon in the linear case.

Indeed, in our context, it is plausible that even a low-dimensional neural network will be able to capture the main features of the time series leading to different best-approximating networks before and after the change.

Because we do not assume that the true auto-regression function is described by this shallow neural network, the errors in the below model can no longer be assumed to be independent but will rather exhibit some time series structure.
More precisely, we replace the above model with 
\begin{equation}\label{eq:quasiModel}
  X_t=
  \begin{cases}
    \Xone_t=f(\Yone_t,\Para_1)+ \e_t^{(1)} & t\le \cp\,,\\
    \Xtwo_t=f(\Ytwo_t,\Para_{2})+ \e_t^{(2)} & t> \cp\,,\\
  \end{cases}
\end{equation}
where $\{\e^{(1)}_t\}$, $\{\e_t^{(2)}\}$ are sequences of stationary and $\alpha$-mixing random variables of polynomial
order and $\Para_1$, $\Para_{2}$ are the best approximating neural network parameters in the sense that $\E(\X_t^{(j)}-f(\Y_t^{(j)},\Para))^2$, $j=1,2$, is minimal (compare \cite{White1989some}). For simplicity,  $\{\Xone_t\}$ and $\{\Xtwo_t\}$ are assumed to be two independent (strictly) stationary
time series, which
differ distributionally, where changes are detectable if the best approximating parameters $\Para_1$ and $\Para_2$ differ.
In an autoregressive case the above assumptions exclude situations where the second time series has starting values from the first time series. While this would certainly be desirable, it greatly complicates proofs and is therefore omitted. While the proof of the  consistency  (Theorem \ref{theo:NLAR-NeuNet_HO_HI} and Corollary \ref{cor:NLAR-NeuNet_CPE_Asymptotics_ConsistentCPE}) only requires an ergodic theorem  which also holds if starting values come from a different distribution, the proof for the rates for the change point estimator (\ref{theo:Op(1)})  require a central limit theorem for the second time series. In case of linear autoregressive models this can easily be achieved by explicitly proving that the influence of the starting values (causing nonstationarity) is asymptotically negligible, for non-linear time series or misspecified situations more complicated arguments such as coupling techniques would be required which are not available in this generality to the best of our knowledge.

In the proofs we will make use of invariance principles as well as \HajekRenyi inequalities, which follow from the following assumptions that have also been used by \cite{KirchTadjuidjeTesting}. However, because we use a universal test statistic capable of detecting all changes in the best approximating parameters, we have to make slightly stronger moment assumptions.
 
\begin{MyAssumption}
\item \label{Ass:MomentAssumptions}
	Let $\{\Xone_t: t \in \Z\}$ and $\{\Xtwo_t: t \in \Z\}$ be stationary time series with $\E\left|\Xone_1\right|^\upsilon < \infty$ and
 $\E\left|\Xtwo_1\right|^\upsilon < \infty$ for some  $\upsilon\geq 3$. Furthermore, assume that $\{\epsilon_t^{(1)}, t\in \Z\}$ and $\{\epsilon_t^{(2)}, t\in\Z\}$ are continuous random variables. 
\item \label{Ass:AlphaMixingAssumptions}
	Let $\{\Xone_t: t \in \mathbb{Z}\}$ and $\{\Xtwo_t: t \in \mathbb{Z}\}$ be independent and  $\alpha$-mixing with
  rate $\alpha(j)=o(j^{-c})$, $c>\upsilon/(\upsilon-2)$ for some $\upsilon>3$.
\item Let the exogenous regressors $\{\V_t: t\in\Z\}$, $\V_t\in \R^d$, be either a deterministic sequence or a stationary time series with finite third moment, $\alpha$-mixing with rate $\alpha(j)$ and independent of the errors  $\{\e_t: t\in\Z\}$ in the true model~\eqref{eq:NLAR-NeuNet_started_auto-regressive_model}.\label{Ass:independent_regpart_stationary}
\end{MyAssumption}

In fact, the latter third moment condition can be relaxed if either the model is correctly specified, or the matrix $A$ in \eqref{eq:NLAR-NeuNet_CPT_definition} guarantees that the test statistic does not involve any partial derivatives with parameters associated to the exogenous regressors. The same holds true for the process.

\subsection{Properties of the neural network estimator}
\label{sec:NLAR_with_NeuNet_Properties_NeuNet}
We approximate the regression function using non-linear
(auto-)regressive models with one-layer neural network as
(auto-)regressive function.

The one-layer neural network $f:\R^{p+d}\times \ParaSpace \rightarrow \R$ (
with $\ParaSpace\subset \R^{(p+d+2)h+1}$ compact) is given by 
\begin{equation}\label{eq:definition_neural_network}
  f(y,\Para)=\nu_0+\sum^h_{i=1}\nu_i\psi(<a_i,y>+b_i)
\end{equation}
with $\nu_0,\,\nu_i,\,b_i\in\R$ and $a_i\in \R^{p+d}$, $i=1,\dots,h$.
\begin{MyAssumption}
\item \label{Ass:differentiableSigmoidFct}
	The activating function $\psi$ is assumed
to be a sigmoid function, i.e. a continuous function fulfilling 
\begin{equation*}
  \psi(y)=1-\psi(-y)\quad \lim_{y\rightarrow \infty}\psi(y)=1\quad \lim_{y\rightarrow -\infty}\psi(y)=0\,.
\end{equation*}
Furthermore, we assume it  is $3$-times differentiable w.r.t. $\Para$ with bounded derivatives, where
  $\Para=(\nu_0,\dots,\nu_h,a_{11},\dots,a_{1(p+d)},a_{21},\dots,a_{h(p+d)},$ $b_1,\dots,b_h)^T$.
\end{MyAssumption}
In this section, we will analyze the least squares estimator for the non-linear \mbox{(auto-)}regressive model \eqref{eq:quasiModel} both under the situation of no change as well as in the presence of a change point. The least squares estimator is defined  as
\begin{equation}
  \label{eq:NLAR-Neural-Network_Estimator}
  \estiPara={\arg\min}_{\Para\in \ParaSpace}Q_n(\Para)
\end{equation}
 with
\begin{equation}\label{eq:NLAR-Neural-Network_least_squares_fct}
  Q_n(\Para):=\sum^n_{t=1}(X_t-f(Y_t,\Para))^2\,,
\end{equation}
for a suitable compact set $\ParaSpace\subset\R^{(p+d+2)h+1}$.

We need the following assumption:

\begin{MyAssumption}
\item \label{Ass:Existence_ParameterMinimizer} 
	For $\lambda$ as in \eqref{eq:NLAR-NeuNet_started_auto-regressive_model} there exists  a  unique minimizer of 
\begin{equation}\label{eq:NLAR-Neural-Network_lossfct}
  \lossFct_{n,\Para}=\lambda\E\left(\Xone_1-f(\Yone_{1},\Para)\right)^2
    +(1-\lambda)\E\left(\Xtwo_1-f(\Ytwo_{1},\Para)\right)^2\,
  \end{equation}
  within  the interior of the
  compact parameter set $\ParaSpace \subset \mathbb{R}^{q}$, with $q=(p+d+2)h+1$.
  \end{MyAssumption}
  Obviously, the uniqueness assumption can only be fulfilled up to the identifiability \footnotetext{Let $\Para=(\nu_0,\mu_1,\dots,\mu_h)$ with $\mu_i=(\nu_i,\alpha_i,\beta_i)$.} 
  conditions for  a non-redundant and irreducible neural network as given by \cite{HwangDing97}, i.e. up to a symmetry transformation\footnote{A symmetry transformation $\pi_k$ of $(\nu_0,\mu_1,\dots,\mu_h)$ is defined as $\pi_k(\nu_0,\mu_1,\dots,\mu_h)=(\nu_0+\mu_k,\mu_1,\dots,\mu_{k-1},-\mu_k,\mu_{k+1},\dots,\mu_h)$. } 
  and transposition\footnote{A transposition (i.e. a permutation changing only two elements)  $\pi$ of $(\nu_0,\mu_1,\dots,\mu_h)$, is given for $i<k$ by $\pi_{ik}(\nu_0,\mu_1,\dots,\mu_h)=(\nu_0,\mu_1,\dots,\mu_{i-1},\mu_k,\mu_{i+1},\dots,\mu_{k-1},\mu_i,\mu_{k+1},\dots,\mu_h)$.}. 
  To guaranty the identifiability, we consider only weight vectors $\Para$ lying in a fundamental domain in the sense of \cite{Rueger1997}. In our case, this is equivalent to $\nu_1\ge \nu_2\ge \dots\ge \nu_h> 0$.

  \begin{MyAssumption}
 \item \label{Ass:Hessmatrix_invertable} Let
\begin{equation*}
M:=\nabla^2 \lossFct_{\limPara}
\end{equation*}
be positive definite.
\end{MyAssumption}

\begin{Proposition}\label{prop:NLAR-NeuNet_Assertion_consistency_normality}
Assume \ref{Ass:MomentAssumptions}, \ref{Ass:independent_regpart_stationary} and \ref{Ass:differentiableSigmoidFct}, \ref{Ass:Existence_ParameterMinimizer}. Then
\begin{enumerate}[label=\alph*)]
\item the parameter estimator is strongly consistent, i.e. 
	\begin{equation} \label{eq:consistentesti}
		\norm{\estiPara-\limPara}\asconv 0 \qquad \text{as }n\rightarrow \infty\,,
	\end{equation}
	where $\limPara$ as in \ref{Ass:Existence_ParameterMinimizer}.
\item  If additionally assumptions \ref{Ass:AlphaMixingAssumptions} and \ref{Ass:Hessmatrix_invertable} hold true, then 
	\begin{equation*}
		\sqrt{n}(\estiPara-\limPara)\dconv \Ncal(0,M^{-1}V M^{-1})\,,
	\end{equation*}
	where 
	\begin{equation*}
		V=\lim_{n\rightarrow\infty}\frac{1}{n}\E\left[\nabla Q_n(\limPara)\transpose{(\nabla Q_n(\limPara))}\right]\,
	\end{equation*}
	and $M$ as in \ref{Ass:Hessmatrix_invertable}.
\end{enumerate}
\end{Proposition}

The second part of the proposition gives in particular the rate of convergence which is used in the proofs below.

\section{Change point procedures}
In this section we will derive asymptotic properties of both the universal test statistics obtained from the neural network approximation as well as the corresponding change point estimators.
\subsection{Asymptotics of the test-statistic}\label{sec:NLAR-NeuNet_CPT_Asymptotics}
We investigate the 
properties of the following test statistics
\begin{align}\label{eq:NLAR-NeuNet_CPT_definition}
& T_n(\eta,\gamma;A)=\max_{1\le k<n}w_{\eta,\gamma}(k/n)\normEstimator{ S(k;\estiPara)},\\
&\text{where } S(k;\Para):=\frac{1}{2}\nabla Q_k(\Para)=\sum^k_{i=1}\nabla f(\Y_i,\Para)(X_i-f(\Y_i,\Para))\nonumber
 \end{align}
 with $A$ fulfilling the following assumption \ref{assump:decision_matrix}, $\Y_t$ given in \ref{notation:NLAR-NeuNet_combined_regression_vector} and $w_{\eta,\gamma}$ as follows.
 \begin{MyNotation}
\item $w_{\eta,\gamma}=\indicator{\eta<s<(1-\eta)}(s(1-s))^{-\gamma}$, with $(\eta,\gamma)\in [0,\frac{1}{2}]^2$. To simplify the notation we write $w_\gamma(s)$ in the case of $\eta=0$, i.e. for  $w_{0,\gamma}(s)$. \label{notation:weight_function}
\end{MyNotation}

	The test statistics in \cite{KirchTadjuidjeTesting} are special cases, where $A$ has only zero elements except for the very first entry 
\begin{equation}\label{eq:A_KirchTadj_Testing}
A=\begin{pmatrix}
1&0&\cdots & 0 \\
0 &0&\cdots &0\\
\vdots & \multicolumn{2}{c}{\ddots}  & \vdots\\
0 & &\cdots& 0
\end{pmatrix}
\end{equation}	
	(compare definition of $\Para$ in \ref{Ass:differentiableSigmoidFct}). For this entry the gradient is equal to one so that the statistic is based on residuals. This is analogous to the first statistics obtained for linear autoregressive time series which were also based on estimated residuals \cite{Kulperger1985,Horvath1993}. The disadvantage is that the corresponding test statistics can only detect a specific class of changes, namely effectively changes that go along with a change in the mean of the corresponding time series. For a positive definite $A$ in the above statistic on the other hand, the statistics is able to detect any change that goes along with a change in the best approximating parameter of the neural network function, in particular in the correctly specified case it will be able to detect all changes. On the other hand, since the number of parameters in this nonparametric setting can be relatively large, changes that are mainly visible in the mean may be detected with lower power only. Therefore, it makes sense to allow also lower rank matrices $A$ in the statistic, increasing power for certain alternatives while at the same time losing power in others. For a detailed discussion of this effect in the multivariate linear setting, we refer to \cite{kirch2015eeg}, similar effects in high-dimensional situations are also discussed in \cite{AstonKirch2018}. Typically, $A$ is chosen such that the limit does not depend on unknowns such as the inverse of $\Gamma$ in Theorem~\ref{theo:NLAR-NeuNet_HO_HI} a) -- or a matrix of lower rank that also uncorrelates the remaining components in a similar manner. Since usually $\Gamma$ is not known and estimators are often not very precise even in for a moderate number of unknown parameters and even more so if the long-run covariance (as in the misspecified model) is needed, it can also make sense to use a different choice and combine this with re-sampling methods to obtain critical values (for a discussion in the multivariate linear case we refer to \cite{kirch2015eeg}).

 \begin{MyAssumption} 
 \item \label{assump:decision_matrix} Let $A$ be a symmetric and positive semi-definite matrix.
 \item \label{assump:H1} For $\limPara$ as in \ref{Ass:Existence_ParameterMinimizer} it holds
 \begin{equation}
 \normEstimator{\E\left[\nabla f(\Yone_t,\limPara)\left(\Xone_t-f(\Yone_t,\limPara)\right)\right]}>0.
 \end{equation}
 By definition of $\limPara$ this is equivalent to 
\begin{equation}
 \normEstimator{\E\left[\nabla f(\Ytwo_t,\limPara)\left(\Xtwo_t-f(\Ytwo_t,\limPara)\right)\right]}>0.
 \end{equation}
 \end{MyAssumption}
For a positive definite matrix $A$ this assumption is always fulfilled in the correctly specified model making all changes detectable. 

\begin{Theorem}\label{theo:NLAR-NeuNet_HO_HI}
   \begin{enumerate}[label=\alph*),ref=\theTheore\ \alph*)]
   \item \label{theo:NLAR-NeuNet_CPT_limit_distri}	Assume \ref{Ass:MomentAssumptions} -- \ref{Ass:Hessmatrix_invertable} and \ref{assump:decision_matrix}. Under \HO the two series 
	\begin{equation*}
	\Gamma_{ij}=\E[Z_{1i}Z_{1j}] + 2\sum^\infty_{l\ge 2}\E[Z_{1i} Z_{lj}]+\sum^\infty_{l\ge 2}\E[Z_{li}Z_{1j}]\,,
	\end{equation*}
	with $Z_{1i}=q_i(t,\limPara)=\nabla_i f(\Y_1,\limPara)(X_1-f(\Y_1,\limPara))$, converge absolutely.\\
	For $\eta=0$ and $\gamma=\frac{1}{2}$ we have 
	\begin{equation*}
	P(\alpha(\log(n))T_n(0,{1}/{2};A)-\beta(\log(n))\le x)\dconv \exp(-2e^{-x}); 
	\end{equation*}
	for $(\eta,\gamma)\in [0,1/2]^2\backslash \{(0,\frac{1}{2})\}$  we have
	\begin{equation*}
		T_n(\eta,\gamma;A)\dconv \sup_{\eta < s< (1-\eta)}\frac{\normEstimator{W(s)-sW(1)}}{(s(1-s))^\gamma}\,,
	\end{equation*}
	where $\{W(t)\}$ is a Wiener process having a covariance matrix $\Gamma=(\Gamma_{ij})_{1\le i,j\le q}$.
\item \label{theo:NLAR-NeuNet_CPT_consistent_Test} Let assumptions \ref{Ass:MomentAssumptions} --
	 \ref{assump:decision_matrix} and 
	if $\eta<\lambda<(1-\eta)$ (with $\cp=\lfloor \lambda n\rfloor$), then under $H_1$ and Assumption \ref{assump:H1} we have for $(\eta,\gamma)\in [0,1/2]^2\backslash (0,1/2)$
	\begin{equation*}
		T_n(\eta,\gamma;A)\pconv \infty\,\qquad 
	\end{equation*}
	and for $(\eta,\gamma)=(0,1/2)$
	\begin{equation*}
	(\log \log n)^{-1/2}T_n(0,1/2;A)\pconv \infty\, .
	\end{equation*}
	
	\item\label{theo:NLAR-NeuNet_CPT_replacing_A_estimator} Replacing the matrix $A$ with a $\sqrt{n}$-consistent estimator does not change the asymptotics under \HO(\ref{theo:NLAR-NeuNet_CPT_limit_distri}) as well as \HI(\ref{theo:NLAR-NeuNet_CPT_consistent_Test}).
   \end{enumerate}
\end{Theorem}

Usually $A$ is chosen in such a way that $ \Gamma^{-\frac{1}{2}}A \Gamma^{-\frac{1}{2}}=D$ is a diagonal matrix with diagonal elements $\delta_i\in \{0,1\}$. Then, the limit is pivotal with $\|W(s)-sW(1)\|_A$ being the Euclidian norm of a $\mbox{rank}(A)$-dimensional  standard Brownian bridge. Typically, $D=\mbox{Id}$ but also other choices of lower rank can be used such as \eqref{eq:A_KirchTadj_Testing} in   \cite{KirchTadjuidjeTesting}, for a detailled discussion in the multivariate linear setting we refer to \cite{kirch2015eeg}. One advantage of such a choice is that only part of the inverse of the long-run covariance  $\Gamma$  needs to be estimated consistenly, which is already 
highly problematic for a moderate number of parameters in particular in the presence of change points. 
On the other hand if the approximation by the neural network is already good enough the corresponding error terms will be almost uncorrelated, in which case it is sufficient to estimate the covariance matrix rather than the long-run covariance matrix. 
By estimating the inverse of the covariance rather than the long-run covariance matrix, we effectively trade a larger estimation error for a small (but asymptotically non-vanishing) model error.

In \cite{KirchTadjuidjeTesting}, it was shown in simulations that this did indeed lead to better small sample results (obviously asymptotically the model error will dominate as seen by the previous theorem). Furthermore, under alternatives a change will typically contaminate the usual empirical covariance estimator. Therefore, in the below simulations and data analyzes, we use the following estimator $\hat{\Gamma}=(\hat{\Gamma}_{ij})_{i,j=1,\dots,h(p+2)+1}$ for $\Gamma=(\Gamma_{ij})_{i,j=1,\dots,h(p+2)+1}$
\begin{align}\label{eq:splittedVarEstimator}
\hat{\Gamma}_{ij}=\frac{1}{n-h(p+2)-1}\Bigg(&\sum^{k_0}_{t=1}(q_i(t,\estiPara^{1}))(q_j(t,\estiPara^{1}))^T+\sum^n_{t=k_0+1}(q_i(t,\estiPara^{2}))(q_j(t,\estiPara^{2}))^T\Bigg)
\end{align}
with $q(t,\Para):=\nabla f(\Y_t,\Para)(X_t-f(\Y_t,\Para))$, $k_0=\argmax_{1\le k<n}\frac{1}{\sqrt{n}}\norm{ S(k;\estiPara)}_{\hat{\Sigma}^{-1}}$ and $\hat{\Sigma}$ is the covariance matrix estimator. 

\subsection{change point estimator}\label{sec:NLAR-NeuNet_CPE}\label{sec:ConsistencyRate}\label{sec:AsymptoticDistri}

In this section, we will analyze the following \linebreak change point estimator for non-linear (auto-)regressive model \eqref{eq:quasiModel} under the assumption there exists a change (\HI, i.e. $\cp=\lfloor \lambda n \rfloor$ with $\lambda\in(0,1)$).\absatz
The weighted multivariate  change point estimator is given as 
\begin{align}\label{eq:NLAR-NeuNet_CPE_definition}
 \cpesti(\eta,\gamma;A)&=\argmax_{1\le k<n}w_{\eta,\gamma}(k/n)\normEstimator{ S(k;\estiPara)}\,,
 \end{align}
where $A$ is symmetric and positive semi-definite and $w_{\eta,\gamma}$ as in \ref{notation:weight_function}. 

	As a first corollary from the above theorem, we can get the consistency of the proposed change point estimator.
\begin{Corollary}\label{cor:NLAR-NeuNet_CPE_Asymptotics_ConsistentCPE}
	Under assumptions
  \ref{Ass:MomentAssumptions}-\ref{Ass:Hessmatrix_invertable}, \ref{assump:H1},
  the  change point estimator $\cpesti(\eta,\gamma)$ \eqref{eq:NLAR-NeuNet_CPE_definition} is a consistent estimator for $\lambda\in (0,1)$ if $0\le \eta<\min(\lambda,1-\lambda)$, i.e.
\begin{equation}
\frac{\cpesti(\eta,\gamma;A)}{n}\pconv \lambda\,.
\end{equation}
In particular, it follows for any $0\le \eta,\alpha<\min(\lambda,1-\lambda)$
\begin{align*}
	P\left( \cpesti(\eta,\gamma;A)=\cpesti(\alpha,\gamma;A) \right)\to 1.
\end{align*}
\end{Corollary}
The corresponding result 
 for $\eta=0$, $\gamma=0$ and 
 $A=(a_{ij})_{1\le i,j\le r}$, $r=(2+p+d)h+1$, where $a_{11}=1$ and $a_{ij}=0$ for $i\cdot j\neq 1$ was proved in \cite{KirchTadjuidjeTesting}. Observe, while the level of the test statistic is indepentent of the choice of the parameter $\eta$ and $\gamma$, under the alternative, the power of the test depends on the chosen $\eta$, $\gamma$ and underlying change point $\tau$. Insides on simulation of the critical quantile for $(\eta,\gamma)\in [0,1/2]^2\backslash (0,1/2)$ we refer to \cite{FRANKE_et_Al_2022}.

We are now ready to state our first result giving the rate of convergence of the change point estimator. In mean change situations it is well known that this rate is optimal. Furthermore, the result shows that for the  change point estimator the actual choice of boundary parameter $\eta$ has no asymptotic influence as long as the true change point lies well within the boundary, i.e. $\eta<\lambda <1-\eta$. 

\begin{Theorem}\label{theo:Op(1)}
Under the assumptions of Corollary \ref{cor:NLAR-NeuNet_CPE_Asymptotics_ConsistentCPE} we have
  \begin{equation*}
    \cpesti(\eta,\gamma;A) -\cp=\OP(1)\,.
  \end{equation*}
\end{Theorem}

The next theorem derives the asymptotic distribution of the estimator which is of interest for two reasons. First, it shows that the above rate cannot be improved, secondly, it helps us to understand the uncertainty connected with the estimator better.
\begin{Theorem}\label{theo:NLAR-NeuNet_CPE_Distri}
Let  the assumptions of Corollary \ref{cor:NLAR-NeuNet_CPE_Asymptotics_ConsistentCPE} hold and $\xi_i^{(z)}=\nabla f(\Yz_i,\limPara)(\Xz_i-f(\Yz_i,\limPara))-\E[(\nabla f(\Yz_1,\limPara)(\Xz_1-f(\Yz_1,\limPara))]$
  for $z=1,2$ be continuous random variables.
Then
	\begin{equation*}
		\cpesti(\eta,\gamma;A)-\cp\dconv \argmax\{W_s-|s| g(s)
		 \normSquared{\D}_A, \ s\in\Z\}
	\end{equation*}
 with
 \begin{align}\label{eq_change_D}\D&=\frac{1}{1-\lambda}\; \E \left[\nabla f(\Yone_t,\limPara)\left(\Xone_t-f(\Yone_t,\limPara)\right)\right]\\
	 &= -\frac{1}{\lambda}\; \E\left[\nabla f(\Ytwo_t,\limPara)\left(\Xtwo_t-f(\Ytwo_t,\limPara)\right)\right],\notag
 \end{align}
 and
	\begin{equation*}
		W_s=\begin{cases}
			0 &,\ s=0\,,\\
			\transpose{\D} A 
			\, \sum^{-1}_{i=s} \xi_i^{(2)} &,\ s<0\,,\\
			\transpose{\D} A 
			\,\sum^{s}_{i=1} \xi_i^{(1)} &,\ s>0\,,\\
			\end{cases}
	\end{equation*}
	 where
	\begin{equation*}
	g(s)=\begin{cases}
	(1-\gamma)(1-\lambda)+\gamma\lambda &,\ s<0\,,\\
	0 &,\ s=0\,,\\
	\gamma(1-\lambda)+(1-\gamma)\lambda &,\ s>0\,.
	\end{cases}
	\end{equation*}
\end{Theorem}
In the general model \eqref{eq:NLAR-NeuNet_started_auto-regressive_model}, $\xi_i^{(z)}$ is a stationary, centered, $\alpha$-mixing time series with $\alpha(j)=o(j^{-c})$.
If model \eqref{eq:quasiModel} is correct, then $$\xi_i^{(z)}=\nabla f(\Yz_i,\limPara) (\e_i)-\E[\nabla f(\Yz_i,\limPara) (\e_i)].$$ For the test statistic as in \cite{KirchTadjuidjeTesting} with A as in \eqref{eq:A_KirchTadj_Testing}, it holds $\xi_i=\e_i$ such that the limit distribution in Theorem \ref{theo:NLAR-NeuNet_CPE_Distri} coincides with the one from the mean change situation as in \cite{AntochHuskova99}. In fact, this result for the mean change is contained as a special case where the neural network has $h=0$ hidden layers. The same holds true in a regression setup with only exogenous variables that are i.i.d.

\section{Simulations and Application}

\subsection{Simulation}
In this section, we examine how well the asymptotic distribution of the change point estimators fits in small sample sizes. We consider both the correctly specified case, where the regression function is truly given by a neural network, and the misspecified case, where the true regression function is only approximated by the neural network. In the controlled setting we show the importance of choosing the proper statistic.\absatz

We compare the results based on the asymptotic distributions of the estimator. To this end we first discuss the method for simulating the asymptotic distribution. Afterwards, we give the results for the correctly speciefied model, followed by the results for the misspecified model.

\subsubsection{Representative of $\tilde \theta$}\label{sec:SimulationsNLAR-NeuNet-representative}
The limiting distribution of the change point estimators is influenced by the unknown parameter $\limPara$. In order to compare the finite-sample distributions with the asymptotic limit distribution, it is necessary to obtain a sufficiently accurate approximation of the neural network parameter for each of our simulation settings. More precisely, the limit depends on the function $f(\cdot,\limPara)$. Consequently, our primary interest lies in ensuring a good fit of this function, particularly in regions where the observations of the regressors are concentrated.

We use that $\widehat{\theta}_N$ converges to $\limPara$ for $N\to\infty$. Since we are interessted in a well enough approximation, run a simulation study to determine a large enough $N$ so that $\widehat{\theta}_N$ can be considered as an approximation of $\limPara$. We need to assure that this estimator is significantly more precise than the estimator for our small sample of length $n$, i.e.\ $N\gg n$. 
In order to decide how large $N$ should be, we consider the variability of the function estimate in the relevant range by considering
\begin{align*}
	V(\hat{\Para}_N):=	\int \mbox{var}\left(f(\hat{\Para}_N,\mathbf{y})-f(\tilde\theta,\mathbf{y})   \right)^2\,g_{\Y_1}(\mathbf{y})\,d\mathbf{y},
\end{align*}
where $g_{\Y_1}$ is the density of $\Y_1$.
This is effectively the variance part of the MISE, which is used as a substitute because we cannot estimate the bias without knowledge of $\limPara$. 
We estimate $V(\hat{\Para}_N)$ by simulating $M$ independent copies of the time series of interest  (correctly specified, misspecified, with or without change at a given location) with length $(N+p)$, i.e.\ $(X^{(r)}_t,\Y^{(r)}_t)$, $t=-p+1,\ldots,N$, $r=1,\ldots,M$. From each of these time series we calculate $\widehat{\theta}_N^{(r)}$, $r=1,\ldots,M$,  and use
\begin{equation*}
	\widehat{V}(\hat{\Para}_N)=	\frac{2} {M(M-1)}\sum^M_{r=1}\sum^M_{s=r+1}\frac{1}{M-2} \sum^{M}_{\substack{i=1\\ i\neq r,\ i\neq s}}\frac{1}{N}\sum_{t=1}^{N}\Big(f\Big({\estiPara}^{r},\Y_t^{(i)}\Big)-f\Big({\estiPara}^{s},\Y_t^{(i)}\Big)\Big)^2\,.
\end{equation*}
This form of the estimator is clearly related to the U-statistics form of the empirical variance.

\subsubsection{Correctly specified model}\label{sec:SimulationsNLAR-NeuNet-correct}~\\
At first, we consider the correctly specified case. Therefore, we choose a simple auto-regression function, a neural network with 1 hidden layer, the model is given by
\begin{equation*}
  X_t=
  \begin{cases}
    0.5+(1+\exp(0.5*(1+0.7 X_{t-1})))^{-1}+\e_t & t\le \cp\,,\\
    \mu+\alpha(1+\exp(0.5*(1+\beta X_{t-1})))^{-1}+\e_t & t> \cp\,
  \end{cases}
\end{equation*}
We consider different changes, where some will effect the mean and one will not. In detail, we define the following models
\begin{align*}
\text{GAR 1:}&\ \mu=0.1,\ \alpha=1,\ \beta=0.7\,,
&\text{GAR 2:}&\ \mu=0.5,\ \alpha=-1,\ \beta=0.7\,,\\
\text{GAR 3:}&\ \mu=0.5,\ \alpha=1,\ \beta=-0.7\,,
&\text{GAR 4:}&\ \mu=0.5,\ \alpha=-1,\ \beta=-0.7\,.
\end{align*}

The analysis on the corresponding test behavior for the test statistic detecting a change in mean can be found in \cite{KirchTadjuidjeTesting}. Instead of only considering a change in mean, we expand the results for $A$ to be chosen as the corresponding derivative. The expected behaviour is a significant increase in the power for GAR 3 while for the other models, the power should be significant lower.

As expected, we observe that the tests have a significant better performance to detect changes as long as the relevant part is detectable. For details on the test study see attachment \ref{sec:SimulationStudies_corretly}. 

While a high power shows the correct decision in determining whether or not a change exist, for the evaluation of the goodness of the estimator a different method is needed. We compare the asymptotic distribution given by \ref{theo:NLAR-NeuNet_CPE_Distri} with the finite sample distribution.  

\begin{figure}[H]
	\centering
	\subfloat[GAR 1: estimated distributions]{
		\includegraphics[scale=0.3]{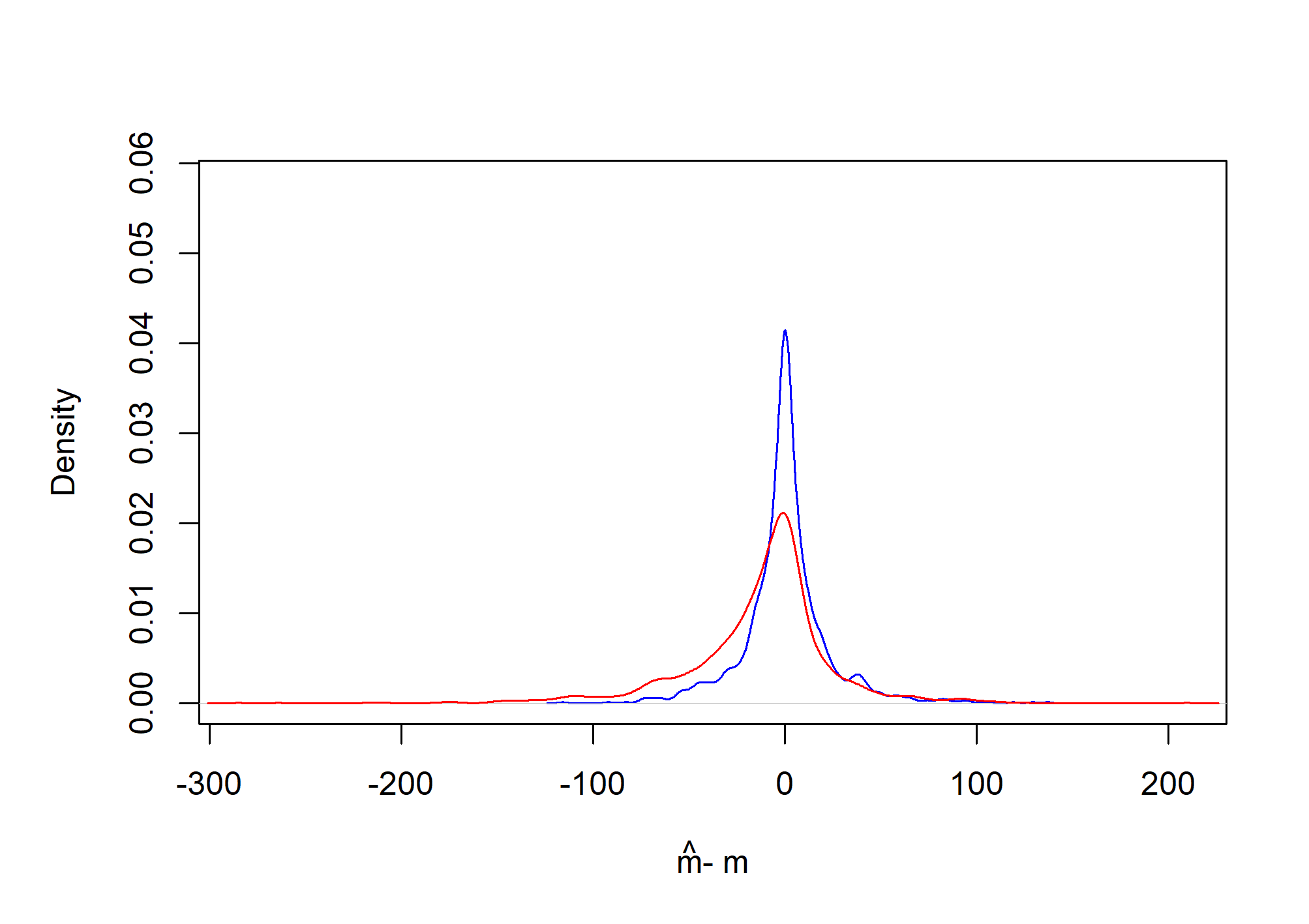}
	}
	\subfloat[GAR 2: estimated distributions]{
		\includegraphics[scale=0.3]{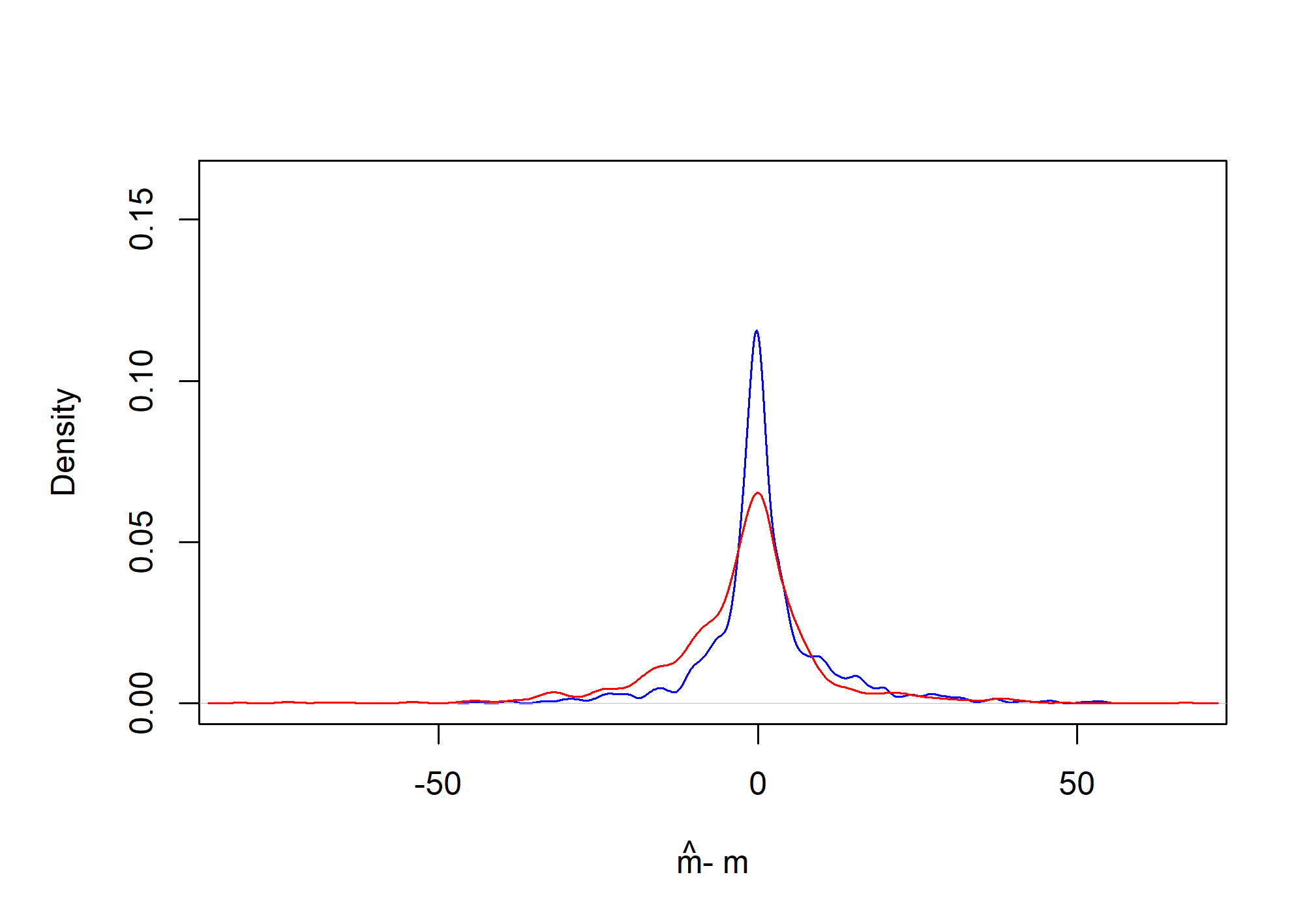}
	}\\
	\subfloat[GAR 3: estimated distributions]{
		\includegraphics[scale=0.3]{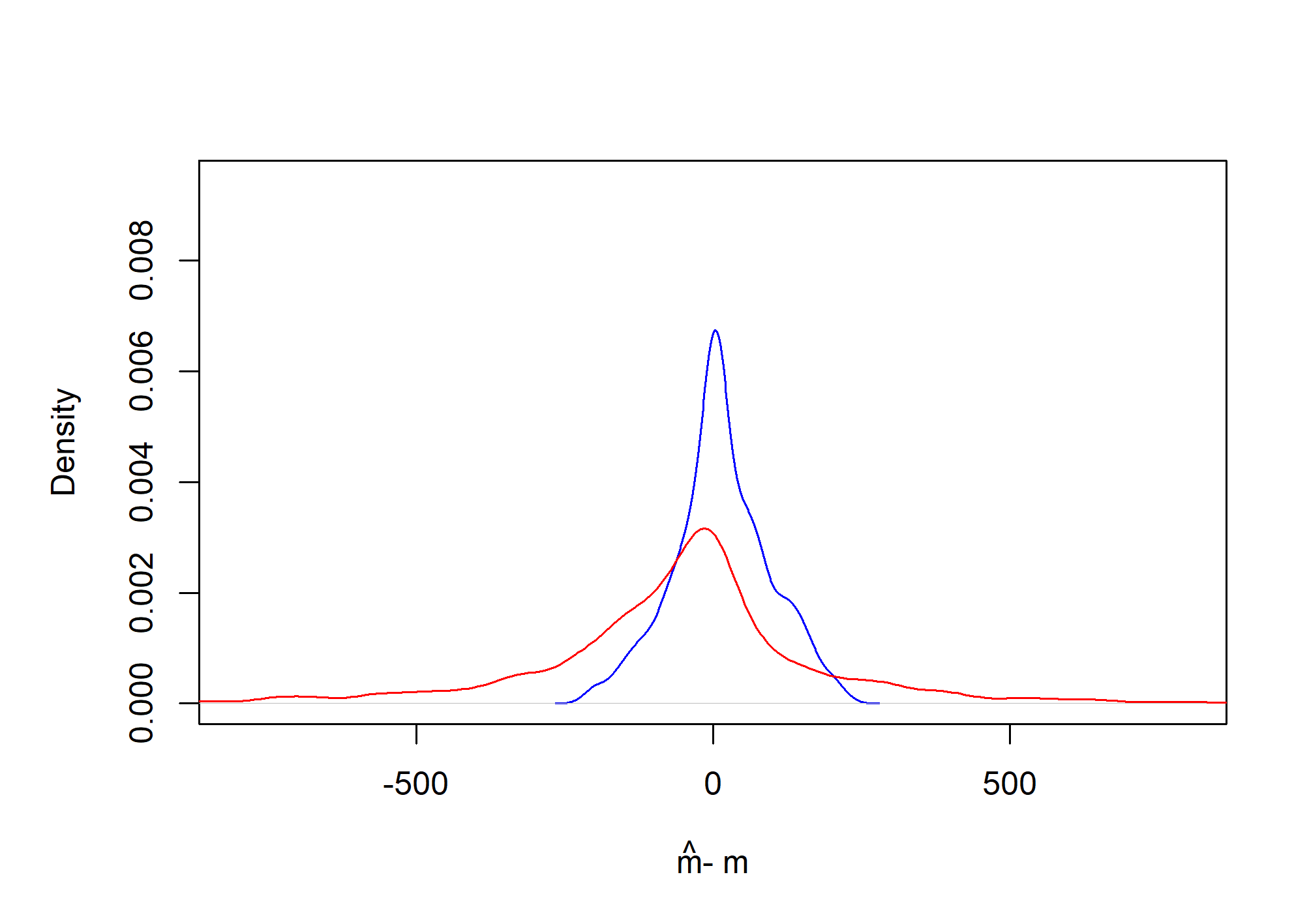}
	}
	\subfloat[GAR 4: estimated distributions]{
		\includegraphics[scale=0.3]{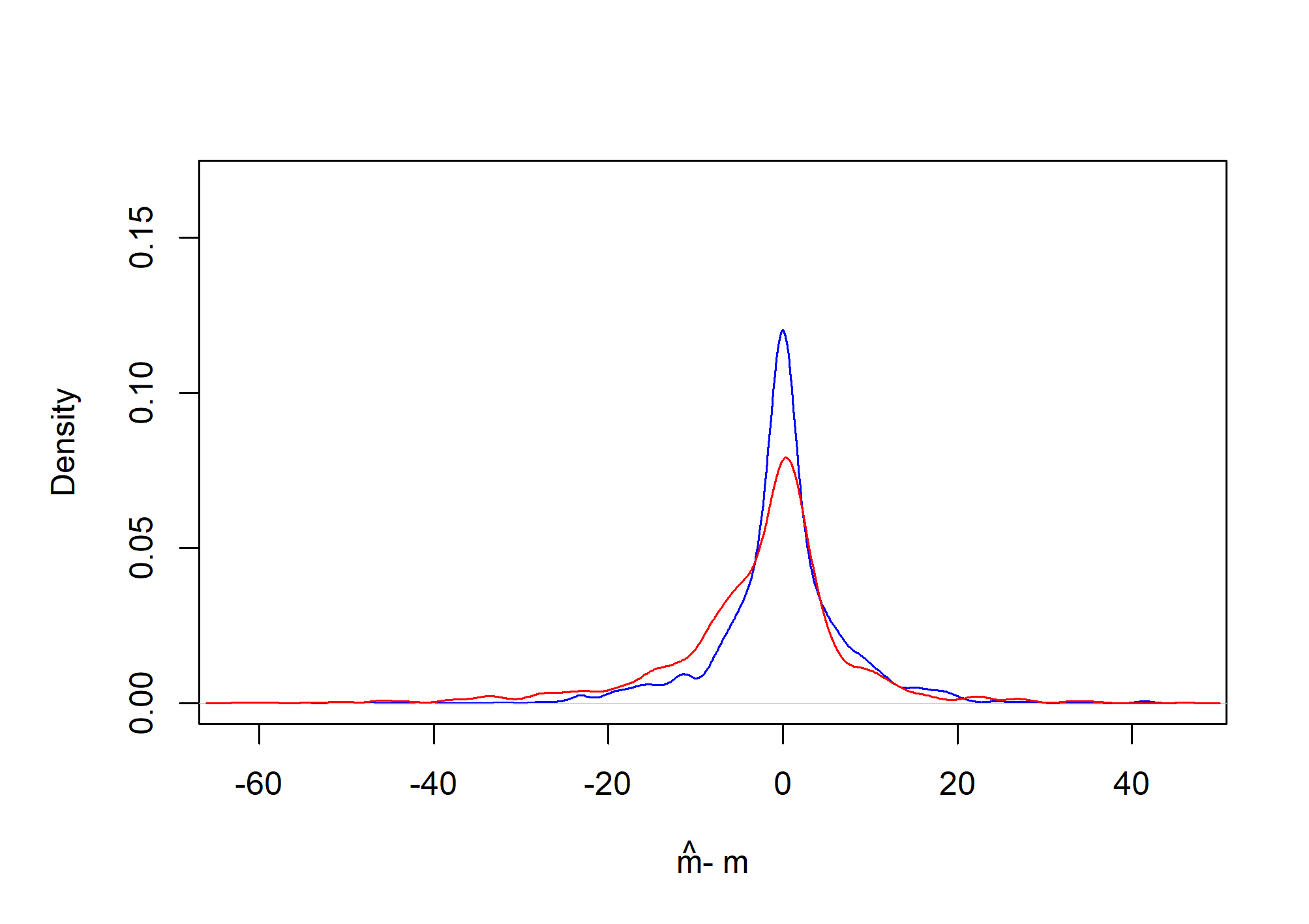}
	}\\[3mm]
	\subfloat[GAR 3 (derivative w.r.t. $\beta$) : estimated distributions]{
		\includegraphics[scale=0.3]{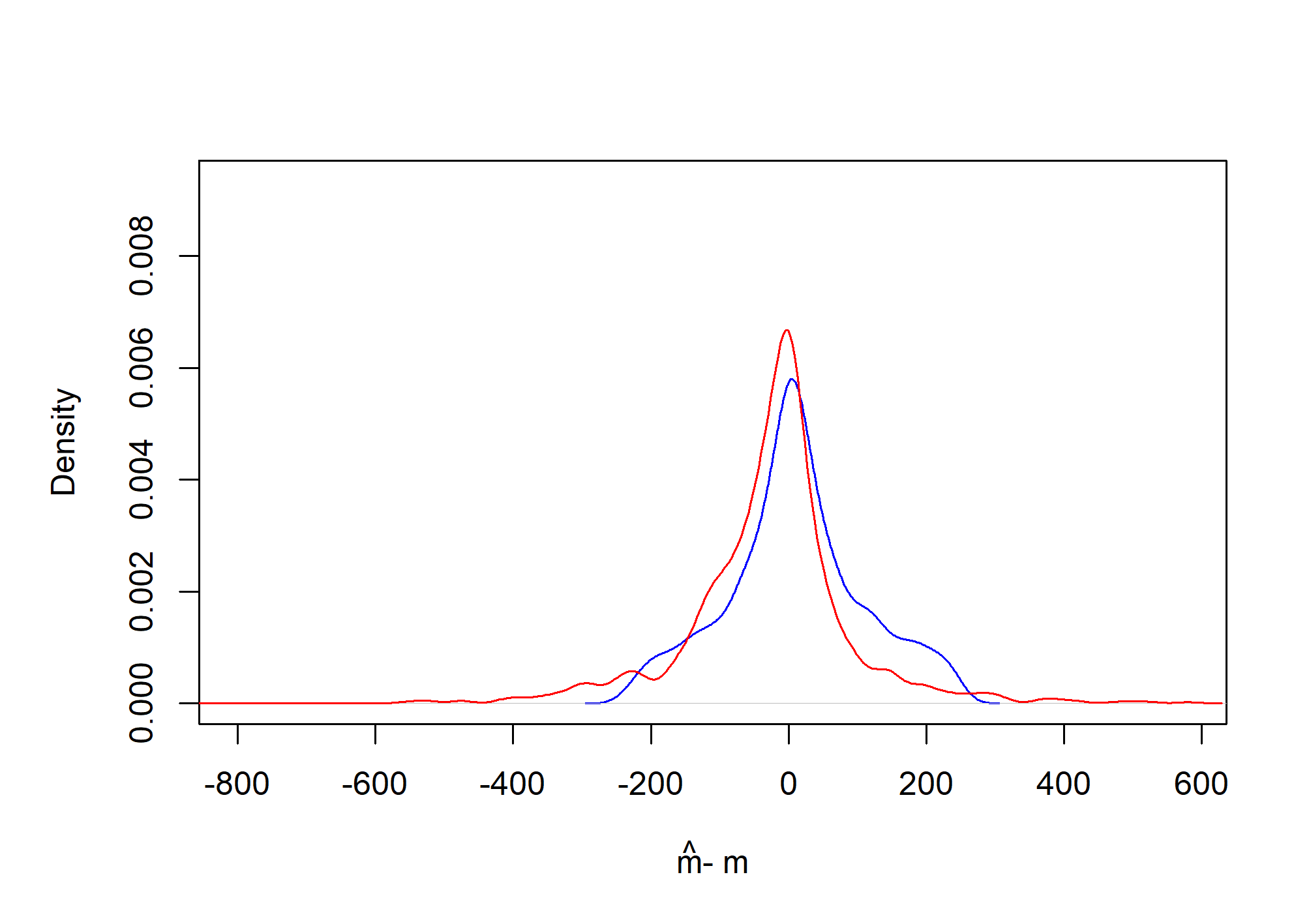}
	}
	\caption{Assumptotic and simulated distribution of the change point estimator for model GAR 1, GAR 2, GAR 3 and GAR 4 with 1 hidden layer neural network as regression function and standard normally distributed errors, N=500}
	\label{fig:Asympt_simulated_distribution_correct_model}
\end{figure}

Figure\ref{fig:Asympt_simulated_distribution_correct_model} shows the comparism of the simulated asymptotic distribution and the finite distribution for $N=500$ and $\tau=0.5$. Even though the tests lag of power for small sample sizes, the asymptotic distributions show for large enough $N$, the estimator will result in reliable estimations.

\subsubsection{Misspecified model}\label{sec:SimulationsNLAR-NeuNet-misspecified}
For the misspecified case, we consider the following situations
\begin{equation*}
  X_t=
  \begin{cases}
    g_1(X_{t-1})+\e_t & t\le \cp\,,\\
    g_2(X_{t-1})+\e_t & t>\cp\,,
  \end{cases}
\end{equation*}
where we choose the autoregression functions as in \cite{KirchTadjuidjeTesting}. The autoregression processes are chosen as
\begin{align*}
\text{AR 1:}&\ g_1(x)=0.3x,\ g_2(x)=0.5+0.1x,
\\
\text{AR 2:}&\ g_1(x)=0.3x,\ g_2(x)=1-0.1x.
\end{align*}
 
For the threshold autoregression models we have 
\begin{align*}
\text{TAR 1:}&\\
\ g_1(x)=&0.3x1_{\{x\ge 0\}}-0.1x1_{\{x<0\}},\
  g_2(x)=(0.5+0.5x)1_{\{x\ge 0\}} -0.3x1_{\{x<0\}},
\\
\text{TAR 2:}&\\
\ g_1(x)=&0.3x1_{\{x\ge 0\}}-0.1x1_{\{x<0\}},\
 g_2(x)=(1-0.1x)1_{\{x\ge 0\}} +(0.5+0.1x)1_{\{x<0\}}.
\end{align*}

The simulation used $N=250,500$ observations and $M=1000$ replications.
In \cite{KirchTadjuidjeTesting} they analyzed the choice of the hidden layer as well as the power of these models. We used the same number of hidden layers for each model. Our calculation, see attachment \ref{sec:SimulationStudies_misspecified}, coincides with the results give in \cite{KirchTadjuidjeTesting}. 

Despite the presence of model misspecification and the known challenges neural networks face in approximating linear functions, our results demonstrate a power close to one, particularly for autoregressive (AR) processes. In the case of mean shifts in AR processes, the test designed to detect changes in the mean shows a distinct advantage. Conversely, for threshold autoregressive (TAR) processes, the derivative-based test exhibits significantly better performance.\absatz
Notably, under correct model specification, where the generalized autoregressive model of order 3 (GAR3) yields a power of less than 50\% for the mean change test statistic, our approach results in substantially higher power, even for TAR processes, regardless of the test statistic employed. These findings suggest that, despite the inherent model misspecification, the neural network-based methods are still capable of delivering robust and reliable results in detecting change points, particularly in the presence of non-linear dependencies.\absatz

\subsection{Application}~\\
As an illustrative example for the application of change point detection, we analyze DAX stock market data from January 2000 to January 2008, a period encompassing the Dot-com bubble. The Dot-com bubble refers to a period of intense speculation on internet-based companies in the late 1990s. After 1995, numerous companies were rapidly founded, often with immediate public stock offerings. High expectations surrounding these companies, particularly those with ".com" domain names, led to significant increases in their stock prices. However, as many of these companies failed to meet speculative expectations, a market correction occurred between 1999 and 2001, during which some companies collapsed entirely. Stock prices continued to decline until 2004, and even after this period, many companies remained undervalued relative to their initial prices.\absatz
In our analysis, we aim to determine whether it is possible to detect and estimate the point of stabilization in the stock market using a non-linear autoregressive (NLAR) model, with change point detection and estimation methods based on neural networks. \absatz

\subsubsection{Fit time series}

\begin{figure}[h]
		\centering
		\subfloat[DAX]{
				\includegraphics[scale=0.3]{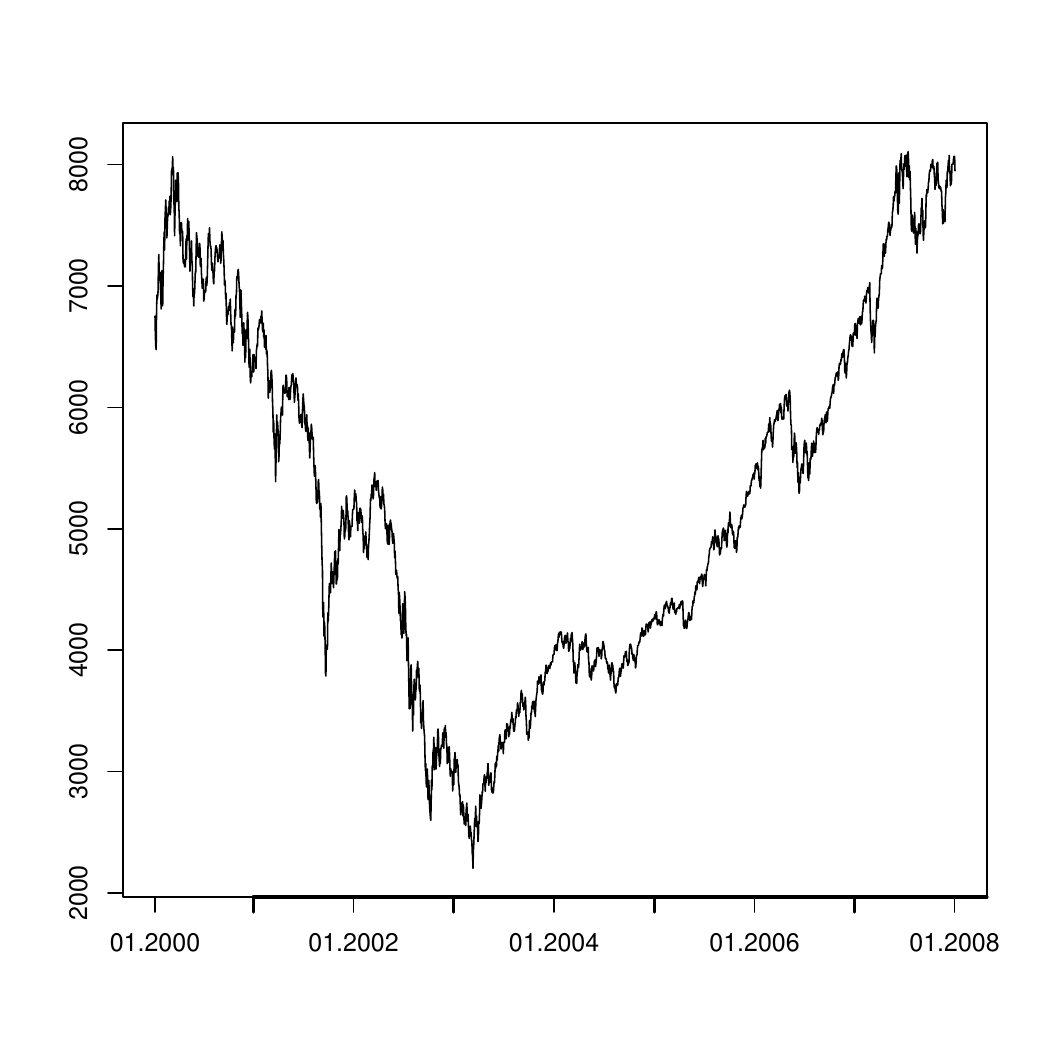}
			}\subfloat[Returns]{
				\includegraphics[scale=0.3]{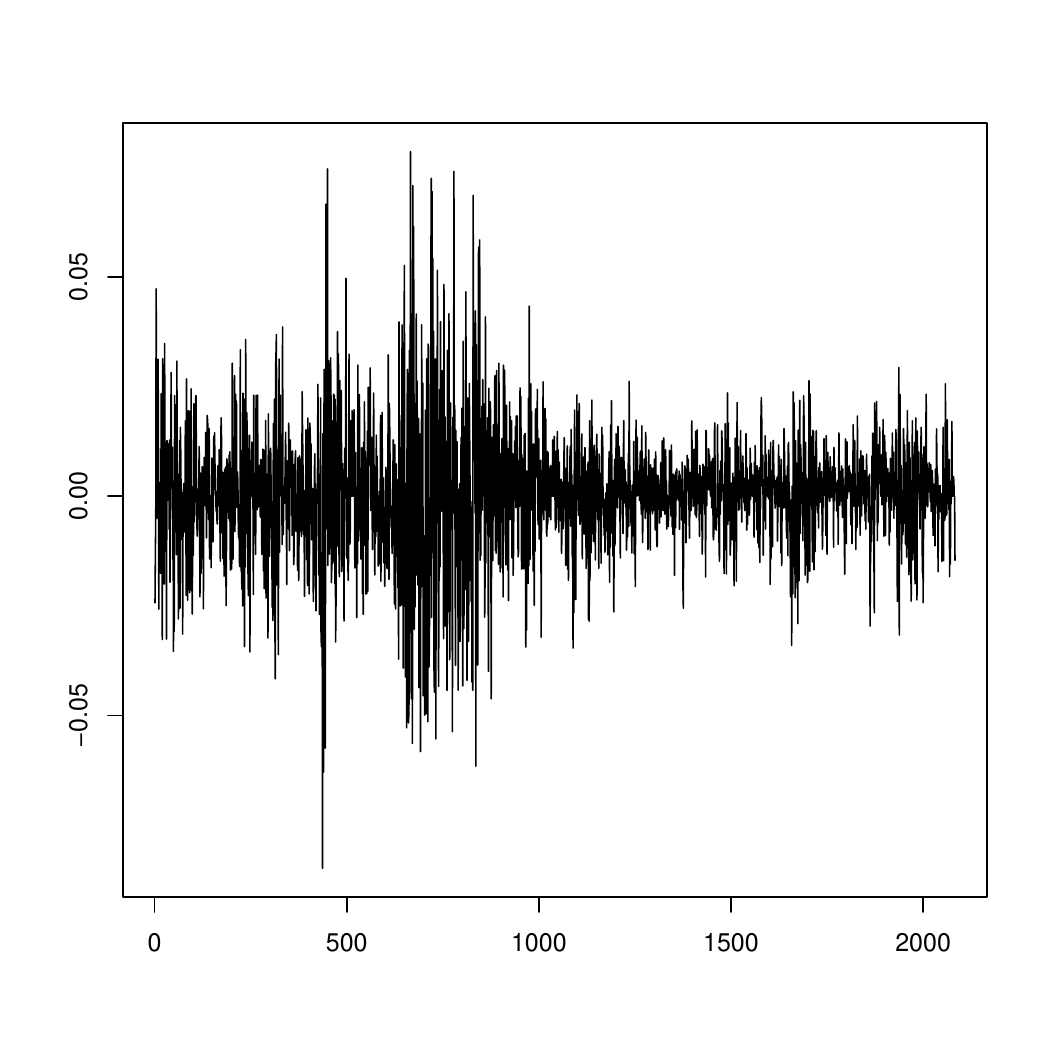}
			}    \caption{DAX-data from 01.2000-01.2008 and the corresponding returns}\label{fig:DAX_data_and_returns}
\end{figure}

To begin the analysis, we fit a non-linear autoregressive (NLAR) model to the DAX data. Let $\{Y_t,t=1,\dots,T\}$ represent the DAX values from January 2000 to January 2008, comprising 2,084 observations. Typically, the focus is on returns, which are calculated as $r_t=(Y_t-Y_{t-1})/Y_{t-1}.$ Figure \ref{fig:DAX_data_and_returns} displays both the DAX values and the corresponding returns for the given time period.

A possible model for the returns is given by $r_t=\sigma_tZ_t$, where $\sigma_t$ is piecewise constant and $Z_t$ follows an autoregressive process, as proposed in \cite{StaricaGranger}. Under this model, the log-returns are expressed as $\log(r^2_t)=2\log(\sigma_t)+\log(Z_t^2)$, implying that a change in variance leads to a change in the mean. However, due to rounding effects, the squared returns can take the value zero. To address this issue, we apply the transformation introduced in \cite{book:Fuller96}, known as the Fuller transformation, which is defined as:
\begin{equation}
	X_t=\log (r_t^2+\rho \hat\sigma^2_r)-\frac {\rho \hat\sigma^2_r}{r^2_t+\rho\hat\sigma^2_r}\,,
\end{equation}
where $\hat{\sigma}^2_r$ is the sample variance of the returns and $\rho>0$  is a weighting parameter. In the following analysis, we use$\rho=0.02$ to obtain a smoother estimate. The Fuller-transformed returns better satisfy the assumptions of the NLAR model.\absatz

For the change point detection (CPT) procedure and the corresponding estimator, a single model order must be selected for the entire time series. It is important to emphasize that order selection methods based on linear autoregressions do not necessarily yield appropriate orders for non-linear autoregressions, such as those modeled by neural networks. Non-linear models typically offer greater flexibility, suggesting that lower orders may suffice. Since our change point approach account for potential model misspecification, we do not require a highly precise model but only a reasonable one. Furthermore, selecting too high an order could result in overfitting, preventing effective detection. To mitigate this risk, we analyze both assumed stationary segments and select the smaller of the two estimated orders.\absatz

We evaluate the first and last 600 observations of the data. Using the \texttt{auto.arima} function from the \textbf{R} package \texttt{forecast}, with the Akaike Information Criterion (AIC) and no moving average component, we fit an autoregressive model to each segment.\absatz

The first segment is best described by an AR(6) model, while the second segment yields an AR(3) model. 
\begin{Verbatim}
ARIMA(6,0,0) with non-zero mean 

Coefficients:
         ar1     ar2     ar3     ar4     ar5     ar6  intercept
      0.0735  0.0575  0.0414  0.0871  0.0757  0.1224    -9.5201
s.e.  0.0405  0.0405  0.0404  0.0406  0.0407  0.0407     0.1371

sigma^2 estimated as 3.394:  log likelihood=-1214.53
AIC=2445.05   AICc=2445.29   BIC=2480.23
\end{Verbatim}
\begin{Verbatim}
ARIMA(3,0,0) with non-zero mean 

Coefficients:
         ar1     ar2     ar3  intercept
      0.0893  0.0417  0.1559   -10.5270
s.e.  0.0404  0.0407  0.0406     0.0981

sigma^2 estimated as 2.967:  log likelihood=-1175.64
AIC=2361.28   AICc=2361.38   BIC=2383.27
\end{Verbatim}
To avoid overfitting, we select the smaller order, $p=3$, for the entire analysis.\absatz
Finally, in the context of neural networks, determining the appropriate number of hidden neurons is crucial. While a detailed discussion of this process is beyond the scope of this paper, further insights can be found in \cite{book:Anders97}. Instead, we experiment with different numbers of hidden neurons in our implementation

\subsubsection{change point detection and estimation}
\begin{figure}[H]
	\centering
	\subfloat[DAX, h=1]{ \includegraphics[scale=0.25]{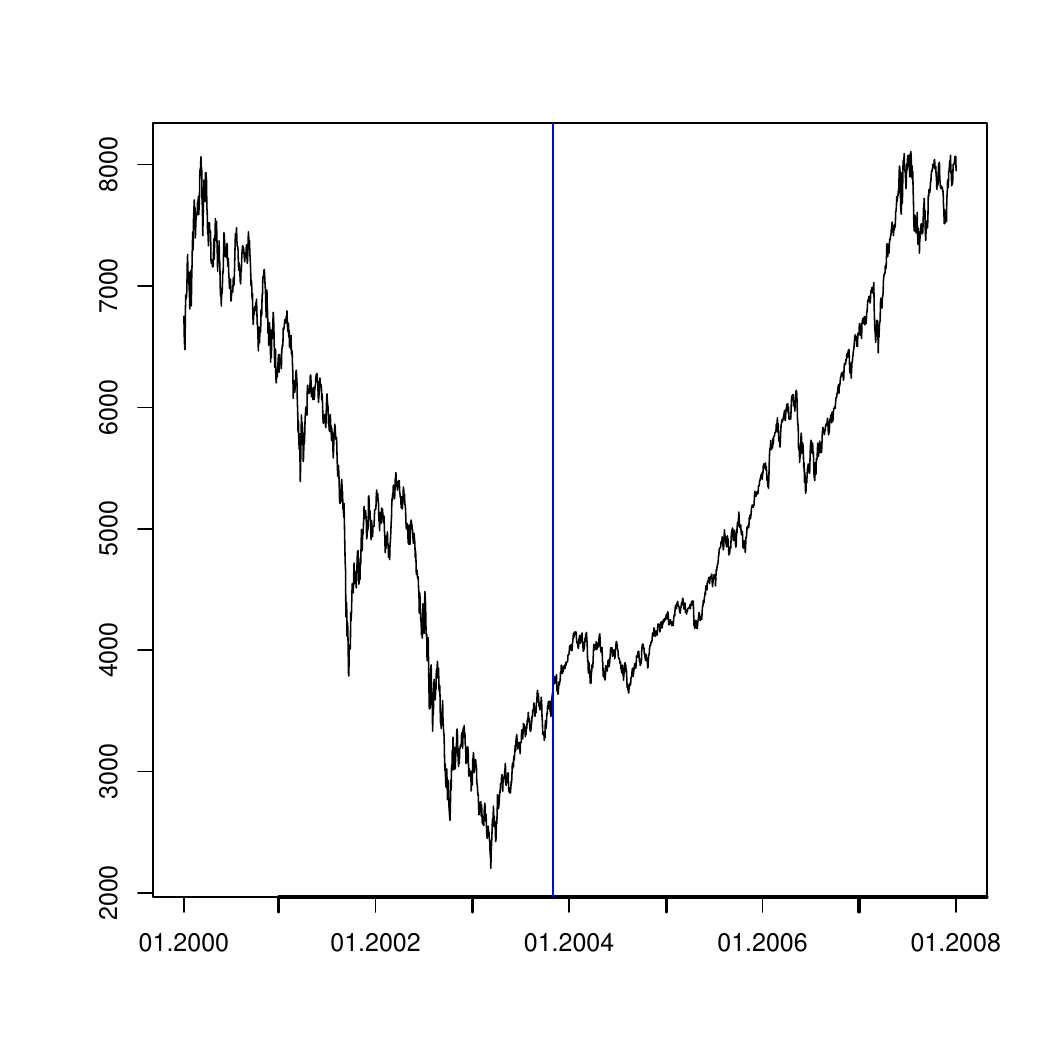}
	}
    \qquad\subfloat[DAX, h=2]{ \includegraphics[scale=0.25]{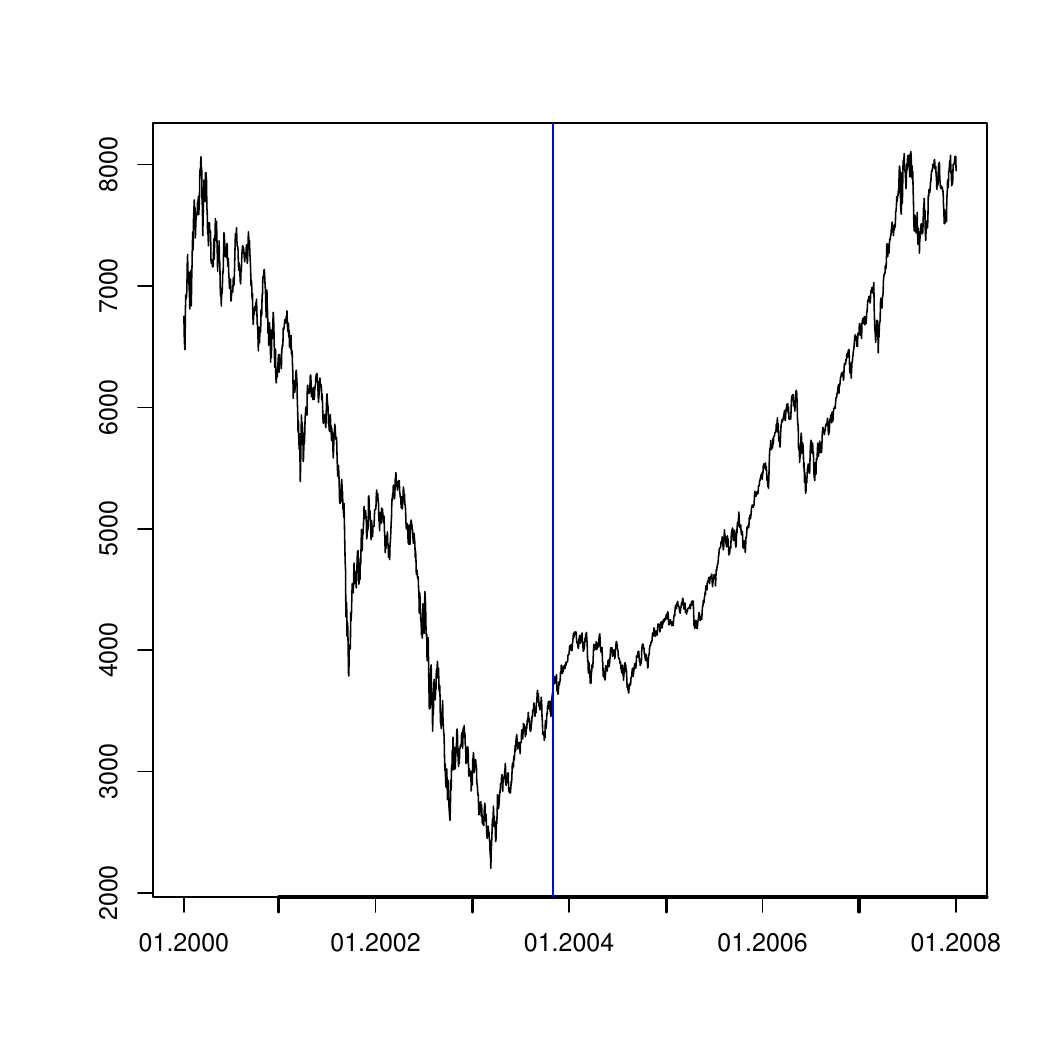}
	}\\
    \subfloat[Fullertransformed returns, h=1]{ \includegraphics[scale=0.25]{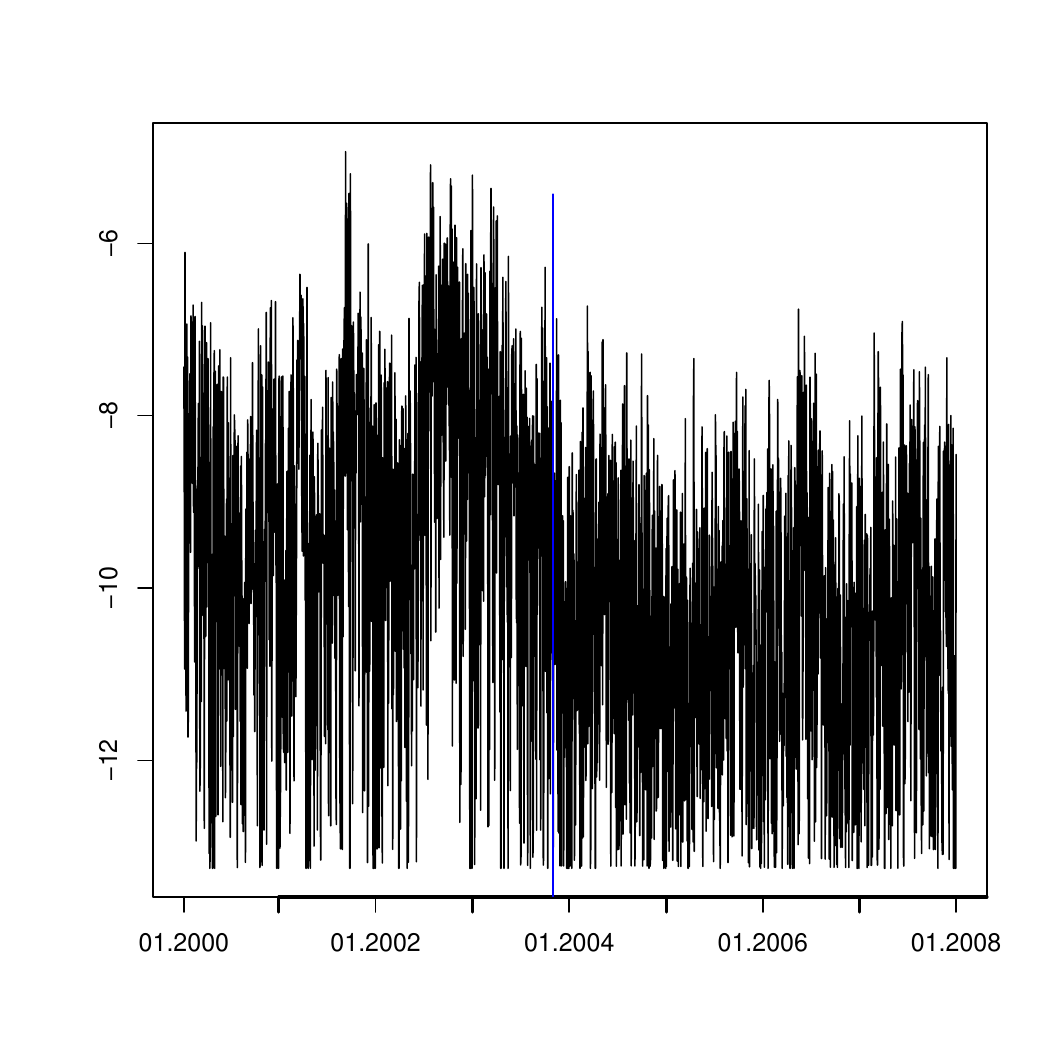}
	}
    \qquad\subfloat[Fullertransformed returns, h=2]{ \includegraphics[scale=0.25]{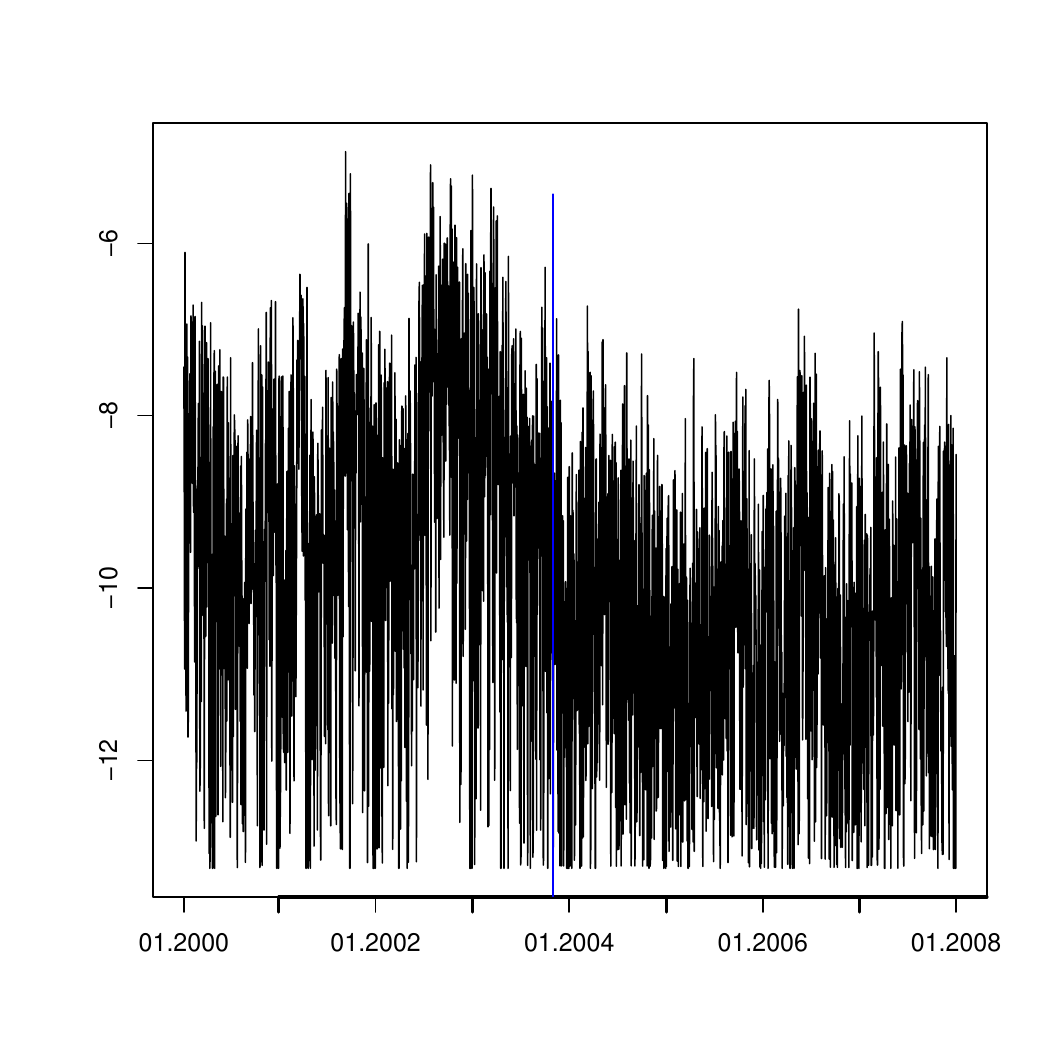}
	}
	\caption{DAX and Fullertransformed returns with the detected and estimated change point (blue line), $\alpha=0.05$}
	\label{fig:Dax_p=3_h=1_2}
\end{figure}

\begin{figure}[H]
	\centering
	\subfloat[DAX, $\alpha=0.01$]{ \includegraphics[scale=0.25]{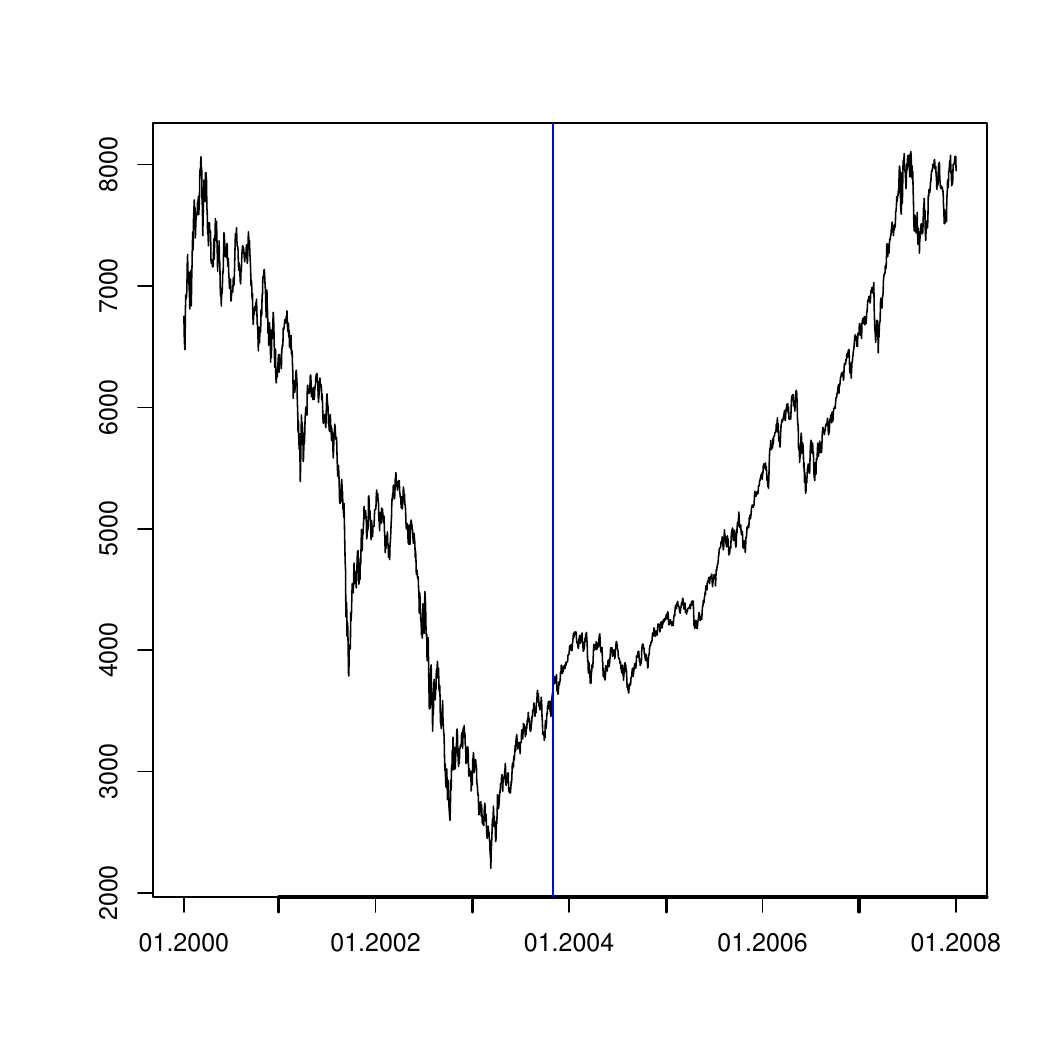}
	}\subfloat[Cumsum, $\alpha=0.01$]{ \includegraphics[scale=0.25]{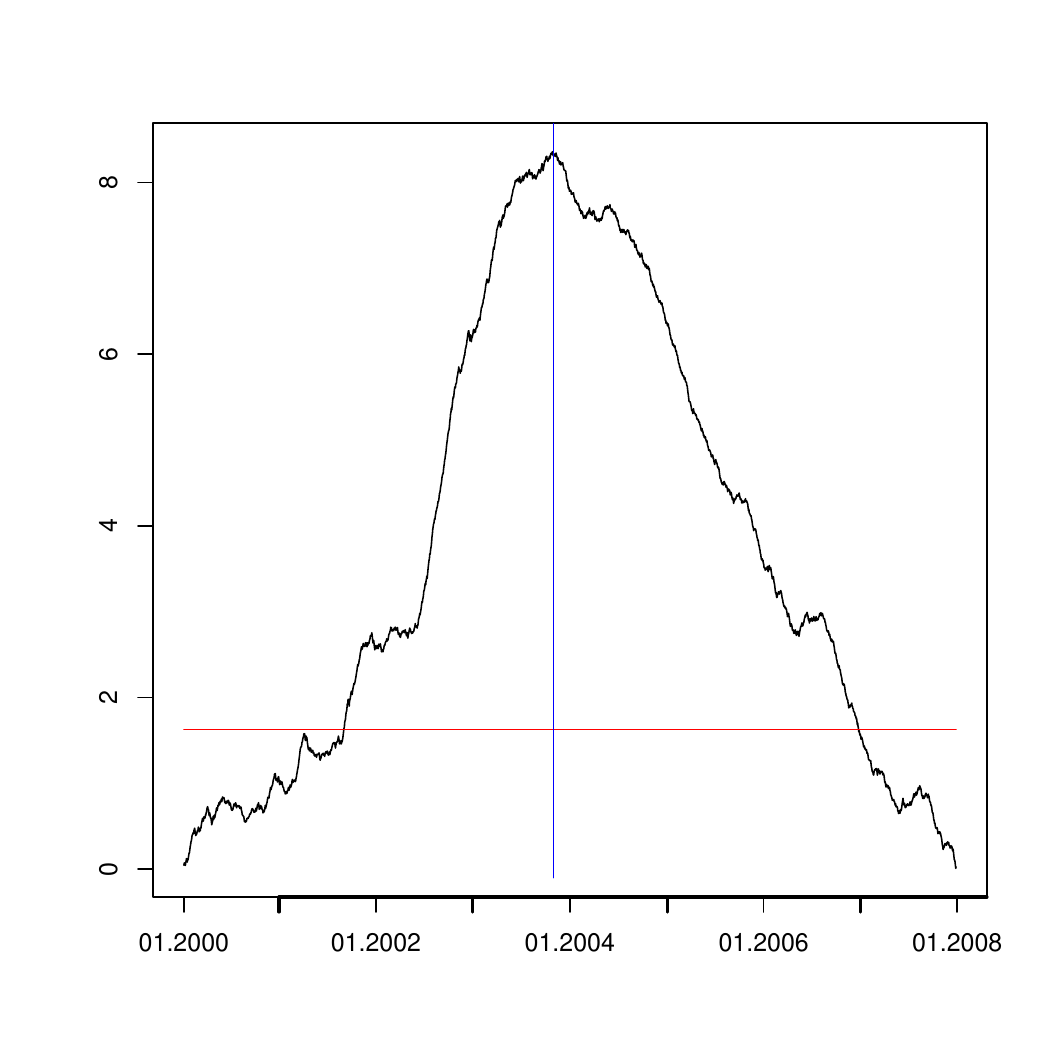}
	}
	\subfloat[returns, $\alpha=0.01$]{ \includegraphics[scale=0.25]{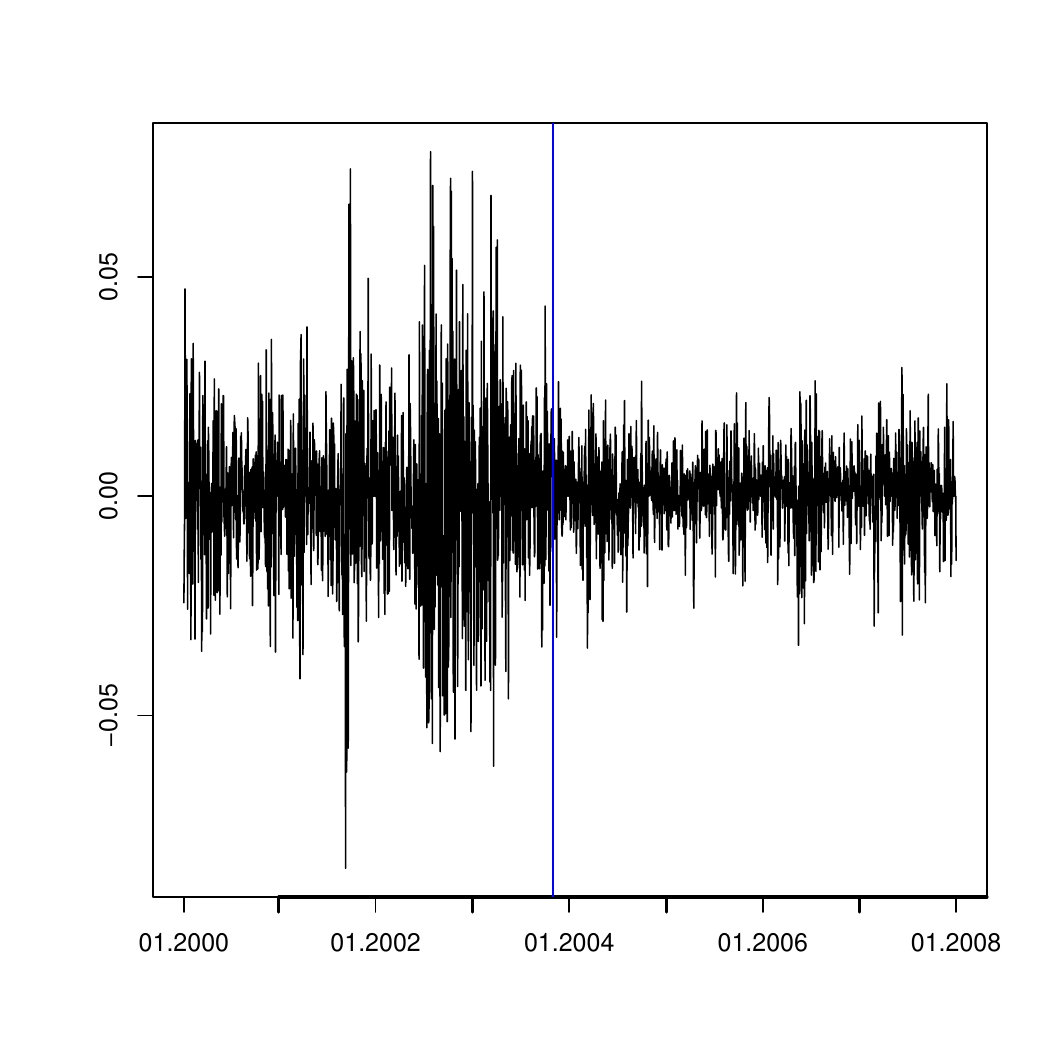}
	}\caption{estimated change point using significance level $\alpha=0.01$ with constant number of hidden neurons $h=1$}
	\label{fig:Dax_diff_sig_levels}
\end{figure}

If we want to use the neural network based change point test and  change point estimator for a NLAR(p)-process, we have to find $p$ and $h$ (number of hidden neurons).
As discussed in the section before, we use a NLAR(3)-model for the test with different numbers of hidden neurons. We choose $h=1,\dots,5$. The results for $h=1,\dots,5$ to the significance level $\alpha=0.05$ do not differ much, we give exemplarily the results for $h=1,2$ in the Figure \ref{fig:Dax_p=3_h=1_2}.\absatz
We mentioned at the beginning, that the stock market for ".com" companies stabilised after 2004. This can be detected at the change in the variance of the returns. To show the significance, we finish this section with the Figure \ref{fig:Dax_diff_sig_levels} giving the DAX, the cumulated sum of the sample residuals (Cumsum, not to be confused with CUSUM, the \teststatistic), the returns and the Fuller transformed returns for the change point test to the significance level $\alpha=0.01$. Even for this level the change point is quite significant, so the dot-com bubble had a significant influence on the stock prices of the DAX.

\section{Conclusion}
In conclusion, the proposed change point test and estimator offers notable advantages for real-world applications. Despite the presence of model misspecification, our approach to change-point detection and estimation delivers consistent and reliable results. The simulation study demonstrates that the derivative-based test provides a significant improvement over the residual-based version, particularly in capturing complex dynamics. Furthermore, the application of neural network-based techniques for both testing and estimation results in plausible and interpretable findings for time change detection.

However, several challenges remain. One key area for future work is the development of strategies to optimize model parameters and tune the neural network for robust and reliable performance without further knowledge about the time series or the change point. Specifically, determining the appropriate number of hidden neurons and the architecture for the network requires further investigation. Another open question involves the optimal choice of the matrix $A$, or devising a systematic approach to determine $A$, which could further enhance the accuracy and efficiency of the procedure. Addressing these issues will not only refine the current methodology but also extend its applicability to a broader range of time series models and real-world datasets.

\section*{Acknowledgment}
We gratefully acknowledge funding from the Ministry of Education and Research Germany (BMBF, grant number 01IS20061 ("EP-KI")) of Stefanie Schwaar. \absatz

\part*{Attachment}
\section{Proofs}\label{sec:NLAR-NeuNet_CPE_Proofs}
\subsection{Proofs of Section 2: Parameter estimator}
\begin{Corollary}
	Under the assumptions \ref{Ass:MomentAssumptions}, \ref{Ass:independent_regpart_stationary} and \ref{Ass:differentiableSigmoidFct} we have
\begin{align*}
	\E[\sup_{\Para\in\ParaSpace}\norm{\nabla f(\Y_t,\Para)(X_t-f(\Y_t,\Para))}]&<\infty\,,&
	\E[\sup_{\Para\in\ParaSpace}\norm{\nabla^2(X_t- f(\Y_t,\Para))^2}]&<\infty\,,\\
	\E[\sup_{\Para\in\ParaSpace}\norm{\nabla^2(\nabla_i f(\Y_t,\Para)(X_t-f(\Y_t,\Para))}]&<\infty\,.
\end{align*}\label{corollary:NLAR-NeuNet_Uniformly_bounded_in_Mean}
\end{Corollary}
\begin{Proof*}
	From the bounded activating function (\ref{Ass:differentiableSigmoidFct}) it follows for all $\theta\in\Theta$ compact
\begin{align*}
	\sup_{\mathbf{y}}\sup_{\theta\in\Theta}|f(\mathbf{y},\theta)|<\infty, \qquad&
\sup_{\theta\in\Theta}|\nabla f(\mathbf{y},\theta)|\leq C \max_{i=1,\dots,d+p}|y_i|,\\
\sup_{\theta\in\Theta}|\nabla^2 f(\mathbf{y},\theta)| < &c \max_{i=1,\dots,d+p}|y_i^2|
\end{align*}
with $C,c>0$ some constants.
The assertion now follows by the moment assumptions \ref{Ass:MomentAssumptions} as well as \ref{Ass:independent_regpart_stationary}.
\end{Proof*}

\begin{Proof}{ of Proposition \ref{prop:NLAR-NeuNet_Assertion_consistency_normality}}
	The proof is analogous to the one of Theorem 3 in \cite{KirchTadjuidjeTesting}. More details can be found in the proof of Theorem 3.2.1 in \cite{Diss:Schwaar}.
\end{Proof}

The following propositions shows uniform convergence of $1/n\,Q_n(\theta)$ and its first and second derivative towards $\lossFct_\Para$.  This result will be needed to prove the null asymptotics as given in Theorem~\ref{theo:NLAR-NeuNet_HO_HI}. 

\begin{Proposition}\label{prop:NLAR-NeuNet_ULN}
Let assumptions \ref{Ass:MomentAssumptions}, \ref{Ass:independent_regpart_stationary} and \ref{Ass:differentiableSigmoidFct} hold and $Q_n(\cdot)$ be as in \eqref{eq:NLAR-Neural-Network_least_squares_fct} and $\lossFct_{\Para}$, $\Para\in \ParaSpace$, as in \eqref{eq:NLAR-Neural-Network_lossfct}. Then for model \eqref{eq:quasiModel} and $\ParaSpace\subset \R^{q}$ ($q=(p+d+2)h+1$) compact, 
we have under \HO as well as under \HI it holds for $n\rightarrow \infty$
\begin{enumerate}[label=\alph*)]
\item \label{prop:NLAR-NeuNet_ULN_parta}
	\begin{equation*}
		\sup_{\Para\in \ParaSpace}\left|\frac{1}{n}Q_n(\Para)-\lossFct_\Para\right|\asconv 0\,,
	\end{equation*}
	\item \label{prop:NLAR-NeuNet_ULN_partb}
	\begin{equation*}
		\sup_{\Para\in \ParaSpace}\norm{\frac{1}{n}\nabla Q_n(\Para)-\nabla \lossFct_\Para}\asconv 0\,,
	\end{equation*}
	with 
	\begin{equation*}
	\nabla \lossFct_\Para=\lambda \E[\nabla f(\Yone_1,\Para)(\Xone_1-f(\Yone_1,\Para))]+(1-\lambda)\E[\nabla f(\Ytwo_1,\Para)(\Xtwo_1-f(\Ytwo_1,\Para))],
	\end{equation*}
	\item  \label{prop:NLAR-NeuNet_ULN_part}
	\begin{equation*}
		\sup_{\Para\in \ParaSpace}\norm{\frac{1}{n}\nabla^2 Q_n(\Para)-\nabla^2 \lossFct_\Para}\asconv 0\,,
	\end{equation*}
	where $\nabla^2$ denotes the Hesse matrix with respect to $\Para$.
\end{enumerate}
\end{Proposition}
\begin{Proof*}
	Due to the assumptions on $\{X_t\}$ and $\{V_t\}$, the proof is analogous to the one of Proposition 1 in \cite{KirchTadjuidjeTesting}, for more details we refer to Proposition 3.2.1 in \cite{Diss:Schwaar}.

\end{Proof*}

\subsection{Proofs of Section 3.1: change point Test}\label{sec:NLAR-NeuNet_CPT_Proofs}
To simplify the notation we introduce the following function.
\begin{MyNotation}
\item Let $q(t,\Para):=\nabla f(\Y_t,\Para)(X_t-f(\Y_t,\Para))$ for $\Para\in \ParaSpace$.
	\label{notation:NLAR-NeuNet_statisticFunction}
\end{MyNotation}
\begin{Proposition}\label{prop:NLAR-NeuNet_CPT_replacement}
	Assume \ref{Ass:MomentAssumptions}- \ref{Ass:Hessmatrix_invertable} and define $q(t,\Para)$ as in \ref{notation:NLAR-NeuNet_statisticFunction}, $\zeta(t)=q(t,\limPara)$. Then under \HO it holds
	\begin{equation*}
		\max_{1\le k<n}\frac{n}{k(n-k)} \norm{\sum^k_{t=1}(q(t,\estiPara)-(\zeta(t)-\sm{\zeta})) 
		}^2=\OP\left(\frac{\log \log n}{n}\right)\,.
	\end{equation*}
\end{Proposition}
The centering with the sample mean $\sm{\zeta}$ takes the additional  variability due to the use of the estimator $\estiPara$ into account. In fact, for $k=n$ the sum over $q(t,\estiPara)$ is $0$ in contrast to the sum over $\zeta_t$. 
\absatz
\begin{Proof*}
	From 
	\begin{align*}
	\norm{\sum^k_{t=1}(q(t,\estiPara)-(\zeta_t-\sm{\zeta}))}^2&=\sum^q_{i=1}\left(\sum^k_{t=1}(q_i(t,\estiPara)-(\zeta_i(t)-\ovl{\zeta}_{i,n}))\right)^2
	\end{align*}
	we get that it is enough to show 
	\begin{equation*}
	\max_{1\le k< n}\sqrt{\frac{n}{k(n-k)}} \left|\sum^k_{t=1}(q_i(t,\estiPara)-(\zeta_i(t)-\sm{\zeta_i})) 
		\right|=\OP\left(\sqrt{\frac{\log \log n}{n}}\right)\,
	\end{equation*}
	for $i=1,\dots,q$. The proof follows analogously to the proof of Lemma 3 in \cite{KirchTadjuidjeTesting} with the  exception that  we take the derivative w.r.t. $\Para$ of the neural network, which is why the higher moment condition is needed. For  a more detailed proof we refer to the proof of Theorem 3.2.2 in \cite{Diss:Schwaar}.
\end{Proof*}
\begin{Proof}{ of Theorem \ref{theo:NLAR-NeuNet_HO_HI}}
	The proof of a) for $(\eta,\gamma)\in[0,1/2]^2\backslash (0,1/2)$ follows from Proposition \ref{prop:NLAR-NeuNet_CPT_replacement} and the invariance principle \cite{KuelbsPhilipp} equivalently as in the proof of Theorem 4 (b) in \cite{KirchTadjuidjeTesting}, see also Theorem 3.2.3 in \cite{Diss:Schwaar}.
	The invariance principle implies a law of iterated logarithm, so the result for $\eta=0$ and $\gamma=1/2$ follows equivalently as in the proof of Theorem 4 (a) in \cite{KirchTadjuidjeTesting}.

	For the proof of b), first consider the case, where $\lambda\geq 0.5$. For $k_0=\lfloor \min(\lambda,1-\eta)\,n\rfloor$ it holds 
		   \begin{align*}
			   S(k_0;\widehat{\theta}_n)=\,\E S(k_0;\limPara)+ [S(k_0;\limPara)-\E S(k_0;\limPara)]+ [S(k_0;\widehat{\theta}_n)-S(k_0;\limPara)]\,,
		   \end{align*}
		   with $S(a;\Para)=\sum^a_{t=1}(q(t,\Para))$. 
The second summand is $\OP(\sqrt{n})$. The third summand is $\oP(n)$, which can be seen by the mean value theorem in combination with a 
uniform law of large numbers \cite[Theorem 6.2]{RangaRao62} and the $\sqrt{n}$-consistency of $\widehat{\theta}$ for $\limPara$ (see Proposition \ref{prop:NLAR-NeuNet_Assertion_consistency_normality}).
		   Hence
	\begin{align}\label{eq:NLAR-NeuNet_CPT_alternative}
			T_n(\eta,\gamma;A)&\geq \frac{1}{\sqrt{n}}\normEstimator{\sum^{k_0}_{t=1}\nabla f(\Y_t,\estiPara)\left(X_t-f(\Y_t,\estiPara)\right)}\,\\\nonumber
		&	= \frac{k_0}{\sqrt{n}}\,\normEstimator{\E\left[\nabla f(\Y_t,\limPara)\left(X_t-f(\Y_t,\limPara)\right)\right]}+\OP(\sqrt{n})\to \infty,
	\end{align}
	concluding the proof.
	For $\lambda<0.5$, the arguments are analogous by chosing $k_0=\lfloor \max(\eta,\lambda) n\rfloor$ and using that by definition of $\widehat{\theta}_n$ it holds $S(k_0;\widehat{\theta}_n)=-\sum_{t=k_0+1}^nq(t,\widehat{\theta}_n)$.
To prove c), first note that by consistency of $\hat{A}_n$ it holds $\|\hat{A}_n^{1/2}A^{-1/2}-\mbox{Id}\|=\oP(1)$. Hence
	\begin{align*}
	\Big|T_n(\eta,\gamma;\hat{A}_n) - &T_n(\eta,\gamma,A)\Big|\\
	&\le \max_{1\le k<n} w(k/n) \left\|\hat{A}_n^{\frac{1}{2}}A^{-1/2} (A^{1/2} S(k;\estiPara)-A^{\frac{1}{2}} S(k;\estiPara)\right\|\\
	&\le \|\hat{A}_n^{1/2}A^{-1/2}-\mbox{Id}\| \, T_n(\eta,\gamma;A)=\oP(1).
	\end{align*}
	The result for the consistency follows because by simular arguments
	\begin{align*}
		\left\|\E\left[\nabla f(\Y_t,\limPara)\left(X_t-f(\Y_t,\limPara)\right)\right]\right\|_{\hat{A}_n}=\left\|\E\left[\nabla f(\Y_t,\limPara)\left(X_t-f(\Y_t,\limPara)\right)\right]\right\|_A+\oP(1),
	\end{align*}
	hence the assertion follows analogously to the proof of b).
\end{Proof}

\subsection{Proofs of Section 3.2: change point Estimator}
\subsubsection{Preliminary results}

\begin{Lemma}\label{lem:momentfromalpha}
  Let $\{Z_t\}$ be $\alpha$-mixing with $\alpha(j)=o(j^{-c})$ for some
  $c>1$.  If there exists $\phi>0$ and $D>0$ such that
  \begin{equation*}
    \E|Z_t|^{2+\phi}\le D\quad \forall t\in\Z\,,
  \end{equation*}
  and
  \begin{equation}
    \alpha(j)=o(j^{-c})\text{ and }
    \phi>\frac {2}{c-1},\label{lemYoko:cond1}
  \end{equation}
  then there exists $\delta\in(0,1]$ and $\Delta>0$ with
  $\phi=\delta+\Delta$ and
  \begin{equation*}
    \E\left|\sum^n_{t=1}Z_t\right|^{2+\delta}\le \Gamma (D,\alpha,\delta,\Delta)n^{\frac {(2+\delta)} 2}\,,
  \end{equation*}
  where $\Gamma$ is a function depending on $D$, $\alpha$, $\delta$
  and $\Delta$.
\end{Lemma}
\begin{Proof*}
The proof follows by Theorem 1 of \cite{Yokoyama}. For details we refer to Lemma C.2.3
in \citep{Diss:Schwaar}.
\end{Proof*}
\begin{Lemma}\label{lem:HajekRenyiAlphaMixing}
  Let $\{Z_t\}$ be a centered, stationary, $\alpha$-mixing process of polynomial order $c>1$
  and assume the ($2+\phi$)-moments exist, with $\phi>\frac 2
  {c-1}$. Then exist $\delta\in (0,1]$ and $\Delta>0$, such that
  $\phi=\delta+\Delta$. Moreover, for any $n\ge 1$, $1\le m\le n$ and
  any positive decreasing
  sequence $b_k$, we have\\
  \begin{enumerate}[label=\alph*)]\itemsep5pt
  \item $\E\left(\max_{k=1,\dots,n}b_k|S_k|\right)^{2+\delta}\le C
    A(\delta) \sum^n_{k=1}b_k^{2+\delta}k^{ \frac{2+\delta}2}$
  \item $\max_{m\le k\le n}b_k
    |S_k|=\OP\left(\left[\sum^n_{k=m}b_k^{2+\delta}k^{
          \frac{2+\delta}2}\right]^{\frac 1 {2+\delta}}\right)$,
  \end{enumerate}
  where $S_k=\sum^k_{t=1}Z_t$, $C$ is a constant and $A(\delta)$ is
  some constant depending on $\delta$, which is as in Lemma \ref{lem:momentfromalpha}.
\end{Lemma}
\begin{Proof*}
The result a) follows by Theorem B.3 in \cite{DissertationKirch}.\\
With analogue arguments b) can be proven, for details see Lemma C.2.4 in \citep{Diss:Schwaar}.
\end{Proof*}
\begin{Lemma}\label{lem:followings}
  Assume 
  \ref{Ass:MomentAssumptions}--\ref{Ass:Hessmatrix_invertable}, then we get for all
  $\theta\in\Theta$ and $z=1$, $2$
  \begin{equation}
    \label{eq:essentiall}
    \begin{split}
    \frac 1 l \sum^l_{t=1}\bigg(&\nabla f(\Y^{(z)}_t,\theta) f(\Y^{(z)}_t,\theta)-\nabla f(\Y^{(z)}_t,\estiPara) f(\Y^{(z)}_t,\hat\theta_n)\bigg)\\
    &=\E\left(\nabla f(\Y^{(z)}_1,\Para)f(\Y^{(z)}_1,\theta)-\nabla f(\Y^{(z)}_1,\limPara)f(\Y^{(z)}_1,\limPara)\right)+\oas(1)
    \end{split}
  \end{equation}for $l\rightarrow \infty$.
\end{Lemma}

\begin{Proof*}
	First notice that it is enough to analyze each entry of the vectors the sum is taken over. Consider the following decomposition of them. We denote the derivative of $f$ w.r.t. $\theta_j$ by $(\nabla f)_j$, $j=1,\dots,{(2+p+d)H+1}$. For $z=1,2$ we have
  \begin{align*}
    \frac 1 l
    \sum^l_{t=1}\bigg((\nabla f)_j&(\Y^{(z)}_t,\theta)f(\Y^{(z)}_t,\theta)-(\nabla f)_j(\Y^{(z)}_t,\estiPara)f(\Y^{(z)}_t,\estiPara)\bigg)\\
    &=\frac 1 l \sum^l_{t=1}(\nabla f)_j(\Y^{(z)}_t,\theta)f(\Y^{(z)}_t,\theta)-\E\left[(\nabla f)_j(\Y^{(z)}_1,\theta)f(\Y_1^{(z)},\theta)\right]\\
    &\qquad +\E\left[(\nabla f)_j(\Y^{(z)}_1,\theta)f(\Y_1^{(z)},\theta)\right]-\E\left[(\nabla f)_j(\Y^{(z)}_1,\limPara)f(\Y_1^{(z)},\limPara)\right]\\
    &\qquad +\E\left[(\nabla f)_j(\Y^{(z)}_t,\limPara)f(\Y_1^{(z)},\limPara)\right]-\frac 1 l \sum^l_{t=1}(\nabla f)_j(\Y^{(z)}_t,\estiPara)f(\Y^{(z)}_t,\estiPara)\\
    &=:A_1+A_2 +A_3\,.
  \end{align*}
  Assumptions \ref{Ass:AlphaMixingAssumptions}, \ref{Ass:independent_regpart_stationary} give $\{\Y_t^{(z)}\}$ is a sequence of $\alpha$-mixing random vectors. From \ref{Ass:AlphaMixingAssumptions} -- \ref{Ass:differentiableSigmoidFct} we conclude that $(\nabla f)_j(\Y^{(z)}_t,\theta)f(\Y^{(z)}_t,\theta)$ is $\alpha$-mixing. Then with the sLLN (strong law of large numbers) and \ref{Ass:MomentAssumptions} it follows $A_1=\oas(1)$. Since the constant on the right hand of \eqref{eq:essentiall} is $A_2$, it is left to 
  show $A_3=\oas(1)$. This follows from 
  \begin{align}\nonumber
    \Bigg|\E\big[(\nabla f)_j&(\Y^{(z)}_t,\limPara)f(\Y_1^{(z)},\limPara)\big]-\frac 1 l \sum^l_{t=1}(\nabla f)_j(\Y^{(z)}_t,\estiPara)f(\Y^{(z)}_t,\estiPara)\Bigg|\\\nonumber
    &\le\left( \sup_{\theta\in\Theta}\left|\frac 1 n \sum^n_{t=1}
        (\nabla f)_j(\Y^{(z)}_t,\Para)f(\Y^{(z)}_t,\theta)-\E[(\nabla f)_j(\Y^{(z)}_1,\Para)f(\Y^{(z)}_1,\theta)]\right|\right)\\\label{eq:Proof_lem_followings}
        &\qquad +\E{\left|(\nabla f)_j(\Y^{(z)}_1,\limPara)f(\Y^{(z)}_1,\limPara)-(\nabla f)_j(\Y^{(z)}_t,\Para)f(\Y^{(z)}_1,\theta)\right|}\Bigg|_{\theta=\estiPara}\,.
  \end{align}
  By Corollary \ref{corollary:NLAR-NeuNet_Uniformly_bounded_in_Mean} and Theorem 6.5 of \cite{RangaRao62}, the first
  term on the right side in \eqref{eq:Proof_lem_followings} is $\oas(1)$.  As we can interchange limit and integration in the second term,
  the consistency of the parameter estimator $\estiPara$ (Proposition \ref{prop:NLAR-NeuNet_Assertion_consistency_normality}) and the continuity of the activating function and its derivatives, finishes the proof.
\end{Proof*}

\begin{Lemma}\label{lem:HelpConsistRate}
  Let $q(t,\Para)$ be as in \ref{notation:NLAR-NeuNet_statisticFunction} and let Assumptions
  \ref{Ass:MomentAssumptions}--\ref{assump:decision_matrix} hold.
  Denote $\tilde S(a,b;\Para)=\sum^b_{t=a}\left(q(t,\Para)-\E q(1,\Para)\right)$ and $\tilde{S}(k;\Para)=\tilde{S}(1,k;\Para)$.
  Then exists $\varsigma\in(0,1)$ such that for fixed
  $\kappa>0$ we have
  \begin{align}
    \label{eq:term1}
  \max_{1\le k\le \cp-\kappa}\left\|\frac 1 {\cp-k}\tilde S (k+1,\cp;\limPara)\right\|_A&=\kappa^{-\frac \varsigma {4+2\varsigma}}\OP(1)\,,\\
 \label{eq:term2}
    \max_{1\le k\le \cp-\kappa}\left\|\frac 1 {\cp-k}\left(\tilde S(k+1,\cp;\estiPara) -\tilde S(k+1,\cp;\limPara)\right)\right\|_A&=\oP(1)\,, \\    
    \label{eq:term3}
	    \max_{1\le k\le \cp}\frac{1}{n}\left\|\tilde S(k;\estiPara)\right\|_A
	    =\oP(1)\,.
\end{align} 
For all fixed but arbitrary constants $C>0$, we have
  \begin{align}
    \max_{0<\cp-k<C}\left\|\tilde S(k+1,\cp;\limPara)\right\|_A=&\OP(1)\,,\label{eq:lemAsymDist1}\\
    \max_{0<\cp-k<C}\left\|\tilde S(k+1,\cp;\estiPara)-\tilde S(k+1,\cp;\limPara)\right\|_A=&\oP(1)\,,\label{eq:lemAsymDist2}\\
    \max_{0\le \cp-k<C}\left\|\tilde S(k;\estiPara)\right\|_A=& \OP\left(\sqrt
      {n}\right)\,.\label{eq:lemAsymDist3}
  \end{align}

  \end{Lemma}

\begin{Proof*}
For every $i=1,\dots,(2+p+d)h+1$ define
  $\zeta_i(t):=q_i(t,\limPara)-\E\left( q_i(1,\limPara) \right)$. Then this is a
  zero mean, stationary and $\alpha$-mixing time series with mixing rate  $\alpha_\zeta(j)=\alpha_X(j-p+1)$ \cite[volume 1, Remark 1.8, Theorem 2.14]{book:Bradley}.  Furthermore,
   for $\upsilon=2+\varsigma+\Delta$, $\varsigma\in(0,1]$ it holds  
   \begin{align}\nonumber
  &  \E[|\zeta_i(t)|^{2+\varsigma+\Delta}]\\\nonumber
  &= O(1)\,\max_{j=1,\dots,p+d}\E[|\Y_{1,j}X_{1}|^{2+\varsigma+\Delta}]\\\nonumber
  &\quad +O(1)\,\max_{j=1,\dots,p+d}\E[|\Y_{1,j}f(\Y_{1},\limPara)|^{2+\varsigma+\Delta}]+\E(q(1,\limPara))^{2+\varsigma+\Delta}\nonumber\\\label{eq:NLAR-NeuNet_CPE_Lemma_rate_deltas}
   & =O(1)\,
  \end{align}
 by assumption \ref{Ass:MomentAssumptions}, \ref{Ass:AlphaMixingAssumptions}, \ref{Ass:independent_regpart_stationary} as well as boundedness of $f(y,\limPara)$ (\ref{Ass:differentiableSigmoidFct}).

Note that any $\alpha$-mixing time series is also $\alpha$-mixing with the same rates in reverse time. From this and 
Lemma \ref{lem:HajekRenyiAlphaMixing} b) with $\varsigma$ as in \eqref{eq:NLAR-NeuNet_CPE_Lemma_rate_deltas} we have 
\begin{align*}\max_{1\le k\le \cp-\kappa}\frac{1}{\cp-k}\left|\sum^{\cp}_{t=k+1}\zeta_i(t)\right|
      =\kappa^{-\frac \varsigma {4+2\varsigma}}\OP(1)\,.
  \end{align*}
By \ref{Ass:MomentAssumptions} and \ref{Ass:independent_regpart_stationary} we have
  $$\sup_{\theta\in \Theta}\|\nabla
  q_i(t,\theta)\|_A\le D\left(\max_{j=1,\dots,p+d}|\Y_{t,j}^2|+\max_{j=1,\dots,p+d}|\Y_{t,j}^2 X_t|\right),$$
  hence by the mean value theorem
  \begin{align}\nonumber
    |q_i(t,\estiPara)-q_i(t,\limPara)|_A&=|(\nabla q_i(t,\theta))^t(\estiPara-\limPara)|\\
    &\le D\, \left(\max_{j=1,\dots,p+d}|\Y_{t,j}^2|+\max_{j=1,\dots,p+d}|\Y_{t,j}^2 X_t|\right)\
    \|\estiPara-\limPara\|\,,\label{eq:lipschitzIneq}
  \end{align}  
  where $\Y_{t,j}$ denotes the $j$th element of $\Y_t$. 
  By \cite[volume 1, Remark 1.8, Theorem 2.14]{book:Bradley} the time series ${\max}_{j=1,\dots
  ,p+d}|\Y_{t,j}^2|$ as well as $\max_{j=1,\dots,p+d}|\Y_{t,j}^2 X_t|$ is $\alpha$-mixing both in forward and backward time.  Using the LLN for $\alpha$-mixing time-series and the $\sqrt{n}$-consistency of $\estiPara$  (\ref{Ass:Existence_ParameterMinimizer}, \ref{Ass:Hessmatrix_invertable}), we get
  \begin{align}
    \max_{1\le k\le \cp-\kappa}&\left\|\frac 1 {\cp-k}\left(S(k+1,\cp;\estiPara) -S(k+1,\cp;\limPara)\right)\right\|_A\notag\\
    &= O(1)\norm{\estiPara-\limPara}\Big(\max_{1\le k\le \cp-\kappa}\frac 1 {{\cp-k}}\sum^{\cp}_{t=k+1}\max_{j=1,\dots ,p+d}|\Y_{t,j}^2|\notag\\
    &\quad + \max_{1\le k\le \cp-\kappa}\frac 1 {{\cp-k}}\sum^{\cp}_{t=k+1}\max_{j=1,\dots ,p+d}|\Y_{t,j}^2X_t|\Big)\notag\\
    &=\OP\left( \frac{1}{\sqrt{n}} \right). \label{eq_neu_ck} \end{align}
  For \eqref{eq:term3} we get 
  \begin{align*}
           \max_{1\le k\le \cp}\left\|
          \tilde S(k;\estiPara)\right\|_A\,
 \le \left( \max_{1\le k\le \cp}\left\|S(k;\hat \theta_n) -S(k;\tilde \theta)\right\|_A +
  \max_{1\le k\le\cp}\left\|\sum^{k}_{t=1}\zeta(t)\right\|_A\right)=\oP(n)  \end{align*}
  by \eqref{eq_neu_ck} and another application of Lemma \ref{lem:HajekRenyiAlphaMixing}. 

  The equality \eqref{eq:lemAsymDist1} follows from the stationarity of $q(t,\limPara)$ and \eqref{eq:lemAsymDist2} similarly after taking \eqref{eq:lipschitzIneq} into account.

 By \eqref{eq:lemAsymDist1} and \eqref{eq:lemAsymDist2} in order to prove \eqref{eq:lemAsymDist3}  it suffices to consider
  \begin{align*}
	  &\max_{0\le \cp-k<C}\left|\sum^{k}_{t=1}(q_i(t,\estiPara)-\E q_i(1,\limPara))\right| \\
	  &\le  \max_{0\le \cp-k<C}\left|\sum^{k}_{t=1}\zeta_i(t)\right|+\max_{0\le m-k<C}\left|\sum^{k}_{t=1}(q_i(t,\estiPara)-q_i(t,\limPara))\right|\\
	  & =\OP(\sqrt{n})\end{align*}
  where the assertion for the first summand follows by Lemma \ref{lem:HajekRenyiAlphaMixing} b) and the assertion for the second summand by the law of large number in addition to the $\sqrt{n}$-consistency of $\estiPara$ for $\limPara$ (see Proposition \ref{prop:NLAR-NeuNet_Assertion_consistency_normality}).
\end{Proof*}

\begin{Lemma}\label{lem:NLAR-NeuNet_CPE_uniformConvergence_of_sum_Needed_for_pruchaPoetscher}
Under \HI (i.e. $\cp=\lfloor \lambda n \rfloor$ with $\lambda\in (0,1]$) and
  \ref{Ass:MomentAssumptions}--\ref{Ass:Hessmatrix_invertable},
  we have
  \begin{align*}
  \sup_{s\in[0,1]}\left|\frac{1}{n} \normEstimator{S(\lfloor sn\rfloor;\estiPara)}-\normEstimator{E_s(\limPara)}\right|=\oP(1)\,,
  \end{align*}
  with $E_s(\Para)=h(s)\,\frac{1}{1-\lambda}\E[\nabla f(\Y_1,\Para)(X_1-f(\Y_1,\Para))]=h(s)\; \D$ with $\D$ as in \eqref{eq_change_D} and  \begin{align*}
  h(s)=\begin{cases}
  (1-\lambda)\,s, &s\le \lambda\,,\\
  \lambda\,(1-s), & s>\lambda\,.
  \end{cases}
  \end{align*}
\end{Lemma}
\begin{Proof*}
By the triangle inequality it follows
\begin{align*}
\sup_{s\in[0,1]}&\left|\frac{1}{n} \normEstimator{S(\lfloor sn\rfloor;\estiPara)}-\normEstimator{E_s(\limPara)}\right|
 \le \sup_{s\in[0,1]}\normEstimator{\frac{1}{n} S(\lfloor sn\rfloor;\estiPara)-E_s(\limPara)}\,\\
 &\le \sup_{s\in(0,1)}\normEstimator{\frac{1}{n} S(\lfloor sn\rfloor;\estiPara)-\frac{1}{n} S(\lfloor sn\rfloor;\limPara)}
 + \sup_{s\in(0,1)}\normEstimator{\frac{1}{n}S(\lfloor sn\rfloor;\limPara)-E_s(\limPara)}\label{eq:NLAR-NeuNet_CPE_Proof_Lemma_limPara_solution_derivative}\\
 &=\oP(1)\,.\nonumber
\end{align*}
where the assertion for the second term follows from the law of large numbers, while the assertion for the first term follows by \eqref{eq:lipschitzIneq}, the law of large numbers and the $\sqrt{n}$-consistency of $\estiPara$ for $\limPara$ (see Propostion \ref{prop:NLAR-NeuNet_Assertion_consistency_normality}). For $s>\lambda $ the same arguments apply after splitting the sum at the change point $m$.
In this case, it holds by the definition of $\limPara$
\begin{align*}
	E_s(\limPara)=\E[\nabla f(\Y_1,\Para)(X_1-f(\Y_1,\Para))]\, \left( s-\frac{\lambda}{1-\lambda}\, (s-\lambda) \right),
\end{align*}
which finishes the proof on noting that $ s-\frac{\lambda}{1-\lambda}\, (s-\lambda)=\frac{\lambda}{1-\lambda} (1-s)$.
\end{Proof*}
\subsubsection{Main results}~\\

\begin{Proof}{ of Corollary \ref{cor:NLAR-NeuNet_CPE_Asymptotics_ConsistentCPE}}
	The first assertion follows from Lemma \ref{lem:NLAR-NeuNet_CPE_uniformConvergence_of_sum_Needed_for_pruchaPoetscher} in combination with Lemma 3.1 of \cite{book:PoetscherPrucha}.

	The second assertion is a simple consequence as for $\eta<\alpha <\min(\lambda,1-\lambda)$ (w.l.o.g) it holds
\begin{align*}
\P(|&\cpesti (\eta,\gamma;A)-\cpesti(\alpha,\gamma;A)| >0)\nonumber\\
= & \P(\mathsmaller{\frac{\cpesti(\eta,\gamma;A)}{n}}-\lambda<\alpha-\lambda)+\P(\mathsmaller{\frac{\cpesti(\eta,\gamma;A)}{n}}-\lambda>1-\lambda -\alpha))\to 0.
\end{align*}
\end{Proof}

Before we can prove the main results, we need the following lemma:
\begin{Lemma}\label{lem_weights}
For any fixed  $\alpha>0$ and $\alpha n\le k< l\le (1-\alpha)n$ it holds
\begin{align*}
	&(a)\qquad \left|w_{\gamma}^2(l/n)-w_{\gamma}^2(k/n)\right|=O\left(\frac{|k-l|}{n}\right),\\
	&(b)\qquad \left|w_{\gamma}^2(l/n)\,l-w_{\gamma}^2(k/n)\,k\right|=O\left(|k-l|\right),\\
	&(c)\qquad 0< c\le  c\, (l-k)\,n \le w_{\gamma}^2(l/n)\,l^2-w_{\gamma}^2(k/n)\,k^2=O\left(n\,|k-l|\right),
\end{align*}
where  $c$ as well as the O-terms only depend on $\alpha$ and $\gamma$.
\end{Lemma}

\begin{Proof*}
	The upper bounds in (a)--(c) follow from the uniform continuity of $w_{\gamma}(\cdot)$ (resp.\ $w_{\gamma}(\cdot)\,\mbox{id}^j$, $j=1,2$) on $[\alpha,1-\alpha]$. 
	The lower bound in (c) follows by the mean value theorem because the derivative of $w_{\gamma}(\cdot)\,\mbox{id}^2$ is bounded away from zero on $[\alpha,1-\alpha]$.
\end{Proof*}

\begin{Proof}{of Theorem \ref{theo:Op(1)}}
  We have to prove that for every $\epsilon>0$ there exists
  $\kappa>0$ such that
  \begin{equation*}
    P(|\cpesti- \cp|>\kappa)=P(\cpesti<\cp-\kappa)+P(\cpesti>\cp+\kappa)\le \epsilon\,.
  \end{equation*}
  For the  event $0<\cpesti<\cp-\kappa$ it holds by  the second assertion in Corollary~\ref{cor:NLAR-NeuNet_CPE_Asymptotics_ConsistentCPE} for any (fixed) $0<\alpha<\min(\lambda,1-\lambda)$
  \begin{equation*}
    P(\cpesti<\cp-\kappa)=P\left(\max_{\alpha n< k <\cp-\kappa}V_k\ge\max_{\cp-\kappa\le k< (1-\alpha)n}V_k\right)+o(1)\,.
  \end{equation*}
  Furthermore, for $k<m$ it holds with $\tilde{S}$ the centered sums as in Lemma~\ref{lem:HelpConsistRate}

  \begin{align*}
    V_k&=\left\langle -(w_{\gamma}(\cp/n)-w_{\gamma}(k/n))S(k;\estiPara)-w_{\gamma}(\cp/n) S(k+1,m;\estiPara),\right.\\
    &\qquad \qquad\left. w_{\gamma}(\cp/n)S(k+1,\cp;\estiPara)+(w_{\gamma}(\cp/n)+w_{\gamma}(k/n)) S(k;\estiPara)\right\rangle_A\\
    &=\Big\langle -(w_{\gamma}(\cp/n)-w_{\gamma}(k/n))\tilde S(k;\estiPara)-w_{\gamma}(\cp/n) \tilde S(k+1,m;\estiPara)\\
    &\qquad \ - (w_{\gamma}(\cp/n)\cp-w_{\gamma}(k/n)k)\,\E q(1,\limPara),\\
    &\qquad  w_{\gamma}(\cp/n)\tilde S(k+1,\cp;\estiPara)+(w_{\gamma}(\cp/n)+w_{\gamma}(k/n)) \tilde S(k;\estiPara)\\
    &\qquad \ +(w_{\gamma}(\cp/n)\cp+w_{\gamma}(k/n)k)\,\E q(1,\limPara)
    \Big\rangle_A\\
    &= \Big\langle B_{1,1}+B_{1,2}+B_{1,3}\,,\,B_{2,1}+B_{2,2}+B_{2,3}\Big\rangle_A\\
    &=B_1+\ldots+B_6, 
    \end{align*}
where
	\begin{align*}
		&\Big\langle B_{1,1}\,,\,B_{2,2}\Big\rangle_A=(w^2_{\gamma}(k/n)-w^2_{\gamma}(\cp/n))\normEstimatorSquared{\tilde S(k,\estiPara)} =:B_1\\
	&\Big\langle B_{1,1}\,,\,B_{2,1}\Big\rangle_A+\Big\langle B_{1,2}\,,\,B_{2,1}\Big\rangle_A+\Big\langle B_{1,2}\,,\,B_{2,2}\Big\rangle_A\\
	&\quad =-w^2_{\gamma}(\cp/n)\Big\langle \tilde S(k,\estiPara)+\tilde S(\cp;\estiPara),\tilde S(k+1,\cp;\estiPara) \Big\rangle_A  =:B_2\\
	&\Big\langle B_{1,1}\,,\,B_{2,3}\Big\rangle_A+\Big\langle B_{1,3}\,,\,B_{2,2}\Big\rangle_A\\
	&\quad=-2(w^2_{\gamma}(\cp/n)\cp - w^2_{\gamma}(k/n)k)\Big\langle \tilde S(k;\estiPara),\E[q(1,\limPara)]\Big\rangle_A
=:B_3\\
	&\Big\langle B_{1,2}\,,\,B_{2,3}\Big\rangle_A+\Big\langle B_{1,3}\,,\,B_{2,1}\Big\rangle_A\\
	&\quad =2 w^2_{\gamma}(\cp/n)\cp\left\langle \tilde S(k+1,\cp;\estiPara)-\tilde S(k+1,\cp;\limPara)\,,\, \E[q(1,\limPara)]\right\rangle_A \\
	&\qquad + 2 w^2_{\gamma}(\cp/n) \cp\left\langle\tilde S(k+1,\cp;\limPara),\,\E[q(1,\limPara)]\right\rangle_A:= B_4+B_5\\
	&\Big\langle B_{1,3}\,,\,B_{2,3}\Big\rangle_A= -(\cp^2 w^2_{\gamma}(\cp/n)-k^2 w^2_{\gamma}(k/n))\normEstimatorSquared{\E[ q(1,\limPara)]}=:B_6
	\end{align*}
	By Lemma \ref{lem:HelpConsistRate}, Lemma~\ref{lem_weights} 
	there exists $\varsigma$
  such that for chosen large enough $\kappa$, it holds
  \begin{align*}
    \max_{\alpha n < k<\cp-\kappa}&\left|\frac{B_1}{n(\cp-k)}\right| 
    \le  \max_{\alpha n < k<\cp-\kappa}\left|\frac{w^2_{\gamma}(k/n)-w^2_{\gamma}(\cp/n)}{\cp-k}\right|\max_{\alpha n < k<\cp-\kappa}\frac{1}{n}\normEstimatorSquared{\tilde S(k;\estiPara)}\\
    &\qquad\qquad\ \; =O(n^{-1})\,\oP(1)=\oP(1),\\
    \max_{\alpha n< k<\cp-\kappa}&\left|\frac{B_2}{n(\cp-k)}\right|\\
    \le &w^2_{\gamma}(\cp/n)\max_{\alpha n< k<\cp-\kappa}\normEstimator{\frac{\tilde S(k+1,\cp;\estiPara)}{\cp-k}} \max_{\alpha n< k<\cp-\kappa}\frac{1}{n} \normEstimator{\tilde S(k;\estiPara)+\tilde S(\cp;\estiPara)}\\
    =&O(1)\oP(1)\,.
    \end{align*}
    Similarly,  we get 
    \begin{align*}
    \max_{\alpha n< k<\cp-\kappa}&\left|\frac{B_3}{n(\cp-k)}\right|\\
    &\le 2 \max_{\alpha n < k<\cp-\kappa}\left|\frac{w^2_{\gamma}(\cp/n)\cp-w^{2}_{\gamma}(k/n)k}{\cp-k}\right| \max_{\alpha n < k<\cp-\kappa}\normEstimator{\frac{1}{n}\tilde S(k;\estiPara)}\normEstimator{\E q(1,\limPara)}\\
    & = O(1) \oP(1)=\oP(1)\,,\\
    \max_{\alpha n < k<\cp-\kappa}&\left|\frac{B_4}{n(\cp-k)}\right|\\
    & \le 2w^2_{\gamma}(\cp/n)\frac{\cp}{n} \max_{\alpha n < k<\cp-\kappa}\normEstimator{\frac{1}{\cp-k}\left(S(k+1,\cp;\estiPara)-S(k+1,\cp;\limPara)\right)}\normEstimator{\E q(1,\limPara)}\\
    &=O(1)\oP(1)=\oP(1)\,,\\
        \max_{\alpha n< k<\cp-\kappa}&\left|\frac{B_5}{n(\cp-k)}\right|\\
	& \le 2w^2_{\gamma}(\cp/n)\frac{\cp}{n} \max_{\alpha n < k<\cp-\kappa}\normEstimator{\frac{1}{\cp-k}\tilde S(k+1,\cp;\limPara)}\normEstimator{\E q(1,\limPara)}\\
	&=\OP(\kappa^{-\frac{\varsigma}{4+2\varsigma}})\,.
  \end{align*}

	By the lower bound in Lemma~\ref{lem_weights} (c) it holds
	\begin{align*}
		\max_{\alpha \le k<m}\frac{n (m-k)}{|B_6|}=O(1)
	\end{align*}
	as well as $B_6<0$ for $k<\cp$ with $\max_{1\le k <\cp-\kappa}B_6\le c$ for an arbitrary constant $c<0$.

 This implies $\max_{1\le k<\cp-\kappa}\left| \frac{B_1+B_2+B_3+B_4}{B_6}\right|=\oP(1)$, hence
 \begin{align*}
    P(\cpesti < &\cp-\kappa)\\
    &\le P\left(\max_{1\le k <\cp-\kappa}V_k\ge V_{\cp}\right)
    = P\left(\max_{1\le k <\cp-\kappa}\left\{B_6\left(1+\oP(1)+\frac{B_5}{B_6}\right)\right\}\ge 0\right)\\
    &\le P\left(\max_{1\le k <\cp-\kappa}\left|\frac{B_5}{B_6}\right|\ge 1-\tau\right)+P\left(1+\oP(1)\le\max_{1\le k <\cp-\kappa}\left|\frac{B_5}{B_6}\right|\le 1-\tau\right)\\
    &\le P\left(\OP(1)\ge (1-\tau)\kappa^{\frac{\varsigma}{4+2\varsigma}}\right)+o(1)\,,
  \end{align*}
  with $0<\tau<1$ arbitrary. This term becomes arbitrarily small for a sufficiently large $\kappa>0$.\\  
  For $\cpesti>\cp+\kappa$  a similar decomposition leads to the conclusion, using that \linebreak$\sum^k_{i=1}(X_i-f(\Y_i,\hat\theta_n))=-\sum^{n}_{i=k+1}(X_i-f(\Y_i,\hat\theta_n))$.
\end{Proof}

\begin{Proof}{ of Theorem \ref{theo:NLAR-NeuNet_CPE_Distri}}
	To simplify the notation we use $\cpesti:=\cpesti(\eta,\gamma;A)$, as $\eta$, $\gamma$ and $A$ are fixed.\absatz
	As in the proof of Theorem 2 in \cite{Antochetal} 
	 it is sufficient to determine the asymptotic distribution of
  \begin{equation*}
    P(\cp-\kappa\le \cpesti\le \cp+x,|\cpesti-\cp|<\kappa)= P\left(\max_{(k-\cp)\in (-\kappa,x]}V_k\ge \max_{(k-\cp)\in (x,\kappa)}V_k\right)\,,
  \end{equation*}
 where $V_k$ is as in the proof of Theorem \ref{theo:Op(1)}. For the proof we have to analyze $V_k$ for $k\in (\cp-\kappa,\cp)$ and for $k\in (\cp,\cp+\kappa)$. In the first case,  let (w.l.o.g.) $x>0$ and $C>x$ be both fixed but arbitrary. By Lemma \ref{lem:HelpConsistRate} and Lemma~\ref{lem_weights} it holds
  \begin{align*}
    \max_{k\in(\cp-C,\cp)}\left|B_1\right|&=\max_{k\in(\cp-C,\cp)}|w^2_{\gamma}(k/n)-w^2_{\gamma}(\cp/n)|\normEstimatorSquared{\tilde S(k;\estiPara)}\\
    &=\OP\left(n^{-1} \right)\OP(n)=\oP(n)\,,\\
    \max_{k\in(\cp-C,\cp)}\left|B_2\right|&\le w^2_{\gamma}(\cp/n)\max_{k\in(\cp-C,\cp)}\normEstimator{ \tilde S(k;\estiPara)+\tilde S(\cp;\estiPara)}\\
	&\qquad\qquad\qquad\cdot\max_{k\in(\cp-C,\cp)}\normEstimator{\tilde S(k+1,\cp;\estiPara) }  \\
	&=\OP(1)\OP\left(\sqrt{n} \right)=\oP\left( n\right)\,,\\
	\max_{k\in(\cp-C,\cp)}\left|B_3\right|&\le 2\max_{k\in(\cp-C,\cp)}|w^2_{\gamma}(\cp/n)\cp - w^2_{\gamma}(k/n)k|\, \max_{k\in(\cp-C,\cp)}\normEstimator{ \tilde S(k;\estiPara)}\\
	&\qquad\qquad\qquad\cdot \normEstimator{\E[q(1,\limPara)]}\\
	&=\OP(1)\OP(\sqrt{n})=\oP(n)\,,\\	
	\max_{k\in(\cp-C,\cp)}\left|B_4\right|&\le 2 w^2_{\gamma}(\cp/n)\cp\max_{k\in(\cp-C,\cp)}\normEstimator{\tilde S(k+1,\cp;\estiPara)-\tilde S(k+1,\cp;\limPara)}\\
	&\qquad\qquad\qquad\cdot \normEstimator{\E[q(1,\limPara)]} \\
	&=\OP(n)\oP(1)=\oP(n)\,.
  \end{align*}
 
  Finally, it holds for any $C>0$ by the differentiability of $w^2_{\gamma}(\cdot)$ and $m=\lfloor \lambda n\rfloor$
  \begin{align*}
	  \sup_{m-C\le k\le m}\frac{m^2w_{\gamma}^2(m/n)-k^2w_{\gamma}^2(k/n)}{n\,(m-k)}=w^2_{\gamma}(\lambda)\,2\,\lambda\left( 1-\gamma\frac{1-2\lambda}{1-\lambda} \right).
  \end{align*}

  We conclude 
  \begin{align*}
	  \max_{l\in ( -\kappa,0 )}&\frac{\,w^{-2}_{\gamma}(\cp/n)}{\cp \left(1-\frac{\cp}{n}\right) }V_{\cp+l}=\max_{l\in ( -\kappa,0 )}\frac{\,w^{-2}_{\gamma}(\cp/n)}{\cp \left(1-\frac{\cp}{n}\right) }\left( B_5+B_6 +\oP(n) \right)\\
		=&\max_{l \in ( -\kappa,0)}\left(
		2\transpose{\D}A \sum^{-1}_{t=l}\tilde \zeta(\cp+t+1)
		-2|l|\left(\gamma \lambda+(1-\gamma)(1-\lambda) \right)\norm{\D}_A^2 +\oP(1)\right)\,,
	\end{align*}
	with $\D$ as in \eqref{eq_change_D}, $\tilde{\zeta}(t)=q(t,\limPara)-\E[q(1,\limPara)]=\nabla f(\Y_t,\limPara)(X_t-f(\Y_t,\limPara)-$\linebreak $\E[\nabla f(\Y_t,\limPara)(X_t-f(\Y_t,\limPara)]$. Define $\xi^{(1)}_{t}\overset{d}{=} \tilde{\zeta}(-t)$, again a stationary $\alpha$-mixing \timeseries, and
	\begin{align*}
	V_{l}^{(1)}=
			 \transpose{\D}A \sum^{-1}_{t=l}\xi^{(1)}_{t}
			-|l|((1-\gamma)(1-\lambda)+\gamma\lambda) \normSquared{\D}_A\,
	\end{align*}
For $k>m$ one arrives similarly at the process
\begin{align*}
   V_{l}^{(2)}=& \transpose{\D}A \sum^{l}_{t=1}\xi^{(2)}_{t}
			-|l|(\gamma(1-\lambda)+(1-\gamma)\lambda) \normSquared{\D}_A\,,
  \end{align*}
  where $\xi^{(2)}(t)\overset{d}{=}\nabla f(\Ytwo_t,\tilde{\theta})(\Xtwo_t-f(\Ytwo_t,\tilde{\theta}))$.
  Defining  \begin{equation*}
    W_l=
    \begin{cases}     
      V_{l}^{(1)},& l<0\\
       0,&l=0\\
      V_{l}^{(2)},& l>0
    \end{cases}\,,
  \end{equation*}
  the proof can now be concluded as in the proof of Theorem 2 in \cite{Antochetal}.
  For every $\epsilon$ exists $n>N_0$ such that for all fixed but arbitrary $\kappa$ it holds
  \begin{align*}
	p(n,\kappa)=P(\cpesti -\cp \leq x, & |\cpesti-\cp|\leq \kappa)\\
  	& = P(\cp-\kappa \leq \cpesti\leq \cp + x, \cpesti-\cp\leq \kappa)\\
  	& = P\left(\max_{(k-\cp)\in[-\kappa,x]}V_k\geq \max_{(k-\cp)\in(x,\kappa]}V_k\right)\\
  	& =  P\left(\frac{\,w^{-2}_{\gamma}(\cp/n)}{\cp \left(1-\frac{\cp}{n}\right) } \max_{(k-\cp)\in(-\kappa,x)}V_k\geq \frac{\,w^{-2}_{\gamma}(\cp/n)}{\cp \left(1-\frac{\cp}{n}\right) }  \max_{(k-\cp)\in(-\kappa,x)}V_k\right)
  \end{align*}
Define
\begin{align*}
	P\left(\max_{l\in[-\kappa,x]} 2 W_l \geq \max_{l\in(-\kappa,x)}2W_l\right) 
	&=P(\argmax_{l\in[-\kappa,\kappa]} W_l\le x)
	& = p(\kappa)
  \end{align*}
We have
\begin{align*}
	\lim_{n\rightarrow\infty}p(n,\kappa)=p(\kappa)
\end{align*}
since $|\frac{\,w^{-2}_{\gamma}(\cp/n)}{\cp \left(1-\frac{\cp}{n}\right) }V_{\cp+l} - 2 W_l|=o_P(1)$ and so by properties of $o_P$ and $O_P$, we have
$P(|\frac{\,w^{-2}_{\gamma}(\cp/n)}{\cp \left(1-\frac{\cp}{n}\right) }V_{\cp+l} - 2 W_l|>C)<\epsilon\,.$
With this uniform convergence and the existence of the limit for $n\rightarrow\infty$, we finally conclude
  \begin{align*}
  \lim_{n\rightarrow\infty} \lim_{\kappa\rightarrow\infty} P(\cpesti -\cp \leq x) & = \lim_{n\rightarrow\infty}\lim_{\kappa\rightarrow\infty} P(\cpesti -\cp \leq, |\cpesti-\cp|\leq \kappa)\\ 
  	& =\lim_{\kappa\rightarrow\infty}\lim_{n\rightarrow\infty} P(\cpesti -\cp \leq, |\cpesti-\cp|\leq \kappa)\\
  	& = P(\argmax W_l \leq x)\,.
  \end{align*}
\end{Proof}

\section{Results Simulation Studies}\label{sec:SimulationStudie}
\subsection{Representative of $\tilde{\theta}$}

We run a simulation with $M=100$ replications and $N=100, 1000, 10 000, 100 000$. For comparison we choose the mispecified autoregressive model AR1 and TAR2 as well as the correct specified model GAR2.

\begin{center}
	\begin{tabular}{l|l|l|l|l}
		 & $N=10^2$ & $N= 10^3$ & $N= 10^4$ & $N= 10^5$ \\\hline\hline
		$\tau=0.25$ & & & &\\
		AR1 &  23.67642 & 0.9504436 & 0.0008279983 & 0.0000984164\\
		TAR2 &  107.2704 &  0.01021047 & 0.0008583396 & 0.0001313354\\
		GAR2 &  4.228445 & 0.009148033 & 0.001412951 & 0.001193053\\[3mm]\hline
		$\tau=0.25$ & & & &\\
	AR1 &   657.735 & 0.01279447 &  0.0008476919 &  0.0001801981\\
	TAR2 &  0.9199325 & 0.01113022 & 0.001043231 &  0.0001878282\\
	GAR2 & 13.03261 & 0.008724946 &  0.001847671 &  0.0008695016\\	
	\end{tabular}
\end{center}

\subsection{Correctly specified model}\label{sec:SimulationStudies_corretly}
Let $T$ denote the test statistic using the derivative w.r.t. the constant $\nu_0$ and $T_a$ the test statistic using the derivative w.r.t. the parameter $a_i$, $i=1,\dots, H$, see \ref{eq:definition_neural_network}. 
~\\[3mm]
\begin{table}[H]
    \centering
    \begin{tabular}{lllllll}
\toprule
	\textbf{N}	& \textbf{Test}& \textbf{GAR 1} & \textbf{GAR 2} & \textbf{GAR 3} & \textbf{GAR 4} \\\midrule
	250 & $T$ & 0.495\ & 0.952 & 0.064 & 0.991 \\
            & $T_a$ & 0.987 & 0.990  & \textbf{0.990} &  0.976 \\\midrule
	500 & $T$ & 0.850 &  1\ \   & 0.151    &  1  \\
            &  $T_a$& 0.992   & 0.989    & 0.995 &  0.991\\\bottomrule
\end{tabular}
    \caption{sample power with processes having a  change point  at $\tau=0.25$}
\end{table}

\begin{table}[H]
    \centering
\begin{tabular}{llllll}
\toprule
	\textbf{N}	& \textbf{Test}& \textbf{GAR 1} & \textbf{GAR 2} & \textbf{GAR 3} & \textbf{GAR 4} \\\midrule
	250 & $T$ &0.770 & 0.996 & 0.156 & 1 \\
            & $T_a$ & 0.973   &  0.984    & \textbf{0.985} &  0.978 \\\midrule
	500 & $T$ & 0.982 &  1   & 0.268    &  1  \\
            &  $T_a$& 0.991   &  0.993   & 0.996 & 0.995\\\bottomrule
\end{tabular}
    \caption{sample power with processes having a  change point  at $\tau=0.5$}
\end{table}
\subsection{Misspecified model}\label{sec:SimulationStudies_misspecified}

~\\
\begin{table}[H]
    \centering
\begin{tabular}{llllll}
\toprule
	N & AR 1 & AR 2 & TAR 1 & TAR 2 \\\midrule
	250 & 0.595& 0.991 & 0.637 & 0.749\\
  & 0.325& 0.308 & 0.358 & 0.268\\
 
	500 & 0.934 & 1 & 0.968 & 0.991\ \\
  & 0.31 & 0.373 & 0.428 & 0.280\ \\
 \bottomrule
\end{tabular}
    \caption{sample power with processes having a  change point  at $\tau=0.25$}
\end{table}

\begin{table}[H]
    \centering
\begin{tabular}{llllll}
\toprule
	N & Test & AR 1 & AR 2 & TAR 1 & TAR 2 \\\midrule
	250 & $T$ & 0.811 & 0.999 & 0.932 & 0.946\\
	& $T_a$& 0.315 & 0.289 & 0.505 & 0.240\\
	500 & $T$ & 0.991 & 1 & 1 & 0.999\\
  &$T_a$ & 0.323 & 0.306 & 0.666 & 0.234\\
 \bottomrule
\end{tabular}
    \caption{sample power with processes having a change point at $\tau=0.5$}
\end{table}

\bibliographystyle{apalike} 
\bibliography{citations}

\end{document}